\documentclass{article}

\usepackage{pdfpages}
\usepackage[margin=1in]{geometry} 
\usepackage{amsmath,amsthm,amssymb}
\usepackage[utf8]{inputenc}
\usepackage{graphicx}
\usepackage{subfigure}

\newcommand{\R}{\mathbb{R}}
\usepackage[hidelinks]{hyperref}
\usepackage{float}
\usepackage{mathtools,cases}
\usepackage{bbm}
\usepackage{breqn}
\usepackage{multirow}
\usepackage{xcolor}
\usepackage{comment}
\usepackage{soul}
\usepackage{bbold}
\usepackage{bm}
\usepackage{ulem}
\usepackage{algorithm}
\usepackage{algorithmic}

\usepackage[numbers]{natbib}

\usepackage{ulem}
\newtheorem{theorem}{Theorem}[section]
\newtheorem{lemma}[theorem]{Lemma}
\newtheorem{definition}[theorem]{Definition}

\newtheorem{remark}[theorem]{Remark}
\newtheorem{example}[theorem]{Example}

\numberwithin{equation}{section}

 \def \d{\mathrm{d}}

\def\P{\mathbb{P}}
\def\E{{\mathbb{E}}}

\renewcommand{\d}{{\rm d}}

\newcommand{\ma}{\color{magenta}}
\newcommand{\bl}{\color{blue}}
\newcommand{\red}{\color{red}}

\usepackage{enumitem}

\usepackage{mathtools}
\usepackage{authblk}
\mathtoolsset{showonlyrefs}
\title{Trading with propagators and constraints: \\ 
applications to optimal execution and battery storage
}

\author[1]{Eduardo Abi Jaber\thanks{eduardo.abi-jaber@polytechnique.edu. The first author is grateful for the financial support from the Chaires FiME-FDD, Financial Risks and Deep Finance \& Statistics at Ecole Polytechnique.}}
\author[2,3]{Nathan De Carvalho\thanks{nathan.decarvalho@engie.com. The second author is grateful for the financial support provided by Engie Global Markets.}}
\author[1]{Huyên Pham\thanks{huyen.pham@polytechnique.edu. This author  is supported by  the BNP-PAR Chair ``Futures of Quantitative Finance", and by FiME, Laboratoire de Finance des March\'es de l'Energie, and the ``Finance and Sustainable Development'' EDF - CACIB Chair  }}
\affil[1]{Ecole Polytechnique, CMAP}
\affil[2]{Université Paris Cité, LPSM}
\affil[3]{Engie Global Markets}

\begin{document}

\maketitle

\begin{abstract}
Motivated by optimal execution with stochastic signals, market impact and constraints in financial markets,  and optimal storage management in commodity markets, we formulate and solve an optimal trading problem with a general propagator model under linear functional inequality constraints. The optimal control is given explicitly in terms of the corresponding Lagrange multipliers and their conditional expectations, as a solution to a linear stochastic Fredholm equation. We propose a stochastic version of the Uzawa algorithm on the dual problem to construct the stochastic Lagrange multipliers numerically via a stochastic projected gradient ascent, combined with a least-squares Monte Carlo regression step to approximate their conditional expectations. We illustrate our findings on two different practical applications with stochastic signals: (i) an optimal execution problem with an exponential or a power law decaying transient impact, with either a `no-shorting' constraint in the presence of a `sell' signal, a `no-buying' constraint in the presence of a `buy' signal or a stochastic `stop-trading' constraint whenever the exogenous price drops below a specified reference level; (ii) a battery storage problem with instantaneous operating costs, seasonal signals and fixed constraints on both the charging power and the load capacity of the battery.
\end{abstract}

\begin{description}
\item[Mathematics Subject Classification (2010):] 49M05, 49M29, 93E20
\item[JEL Classification:] C61, G10, G14
\item[Keywords:] Stochastic Convex Programming, Constraints, KKT optimality conditions, Lagrange multipliers, Fredholm, Stochastic Uzawa, Battery Energy Storage System, Optimal trading, Market impact and signal
\end{description}



\section{Introduction}

\textit{Constraints} play a pivotal role in real-life optimization problems. Whether originating from regulatory requirements or physical restrictions, constraints significantly impact the decision-making process imposing limitations on the feasible solutions, and in some cases, may even be stochastic. 

In this paper, we formulate and solve, both theoretically and numerically, a class of stochastic control problems with path-dependencies and almost sure constraints. Such class of problems is motivated by two practically relevant situations in quantitative finance and commodity markets where constraints must be taken into account in an uncertain environment.

\begin{itemize}
    \item \textbf{Optimal execution and trading under constraints in traditional financial markets.} 
\textit{Optimal execution} is the art of liquidating an initial long position of a given asset before a prescribed time horizon, while taking into account the trade-off between price signal, volatility and market impact.
Initiated by the seminal work of Almgren and Chriss \cite{almgren2001optimal}, the problem of optimal execution has had many refinements, see for exemple \cite{cartea2015algorithmic} or \cite{gueant2016financial}. The addition of constraints is crucial in practice to prevent excessive trading that could destabilize the market or to comply with regulations. They can vary widely, including constraints on the trading rate, such as restrictions on buying when liquidating, or constraints on the running inventory that prevent short selling. Additionally, index tracking constraints are often applied to ensure that the trading strategy aligns with the performance of a benchmark index, maintaining desired tracking accuracy.  Terminal constraints may require liquidating a specific quantity by a set time or achieving a particular portfolio value, which could be deterministic or stochastic \cite{nutz2023unwinding}, with the possibility of relaxing the full inventory constraint in the case the asset's price moves too adversely, see, for example, \cite{ackermann2024reducing} or \cite{aksu2023optimal}.
 These constraints ensure that execution strategies are both effective and compliant in complex market environments.

\item 
{\textbf{Energy arbitrage \& Battery Energy Storage Systems (BESS) optimization.}
\textit{Energy arbitrage} is the art of buying low, storing and selling high a commodity based on some stochastic price signals in order to maximize expected revenues. One of the key challenges in such operational control problems is to account for both the uncertainty of commodity prices as well as some physical storage constraints and influx/outflux flow restrictions, known in the literature as the \textit{warehouse problem}, see \cite{secomandi2010optimal} and the references therein. The case of Markovian diffusion dynamics for the commodity price, and when neglecting instantaneous operating costs and market impact effects, has already been extensively covered in the literature \cite{carmona2010valuation}, \cite{lai2010approximate}, \cite{secomandi2010optimal}, \cite{jiang2015optimal}, \cite{bachouch2022deep}, \cite{lemaire2024swing}, where typically optimal `bang-bang' policies are derived using a dynamic programming formulation. In practice however, (i) `bang-bang' policies may not be desirable for the physical asset's wear (e.g.~inertia, over-heating), (ii) market impact effects may become non-negligible as the volume of storage increases, and (iii) more realistic stochastic path-dependent price signals can better capture arbitrage opportunities. As a typical example, the optimization of Battery Energy Storage Systems consists in determining the best strategy to manage (charge and discharge) electricity storage systems in order to maximize efficiency, lifespan, and economic benefits\footnote{For example, flexible battery storage of electricity can contribute to better match a seasonal and yet uncertain demand and compensate for unexpected production incidents more effectively by arbitraging the sudden spikes of prices observed on short-term power markets which cannot be captured by mainstream technologies available in the production mix due to a lack of operational flexibility, see \cite{paatero2005effect}.}. Adding to the complex reality of such storage optimization problem, battery assets can be optimized on different markets with various time horizons, including, for example, in Europe, the ancillaries, the day-ahead, the intra-day or the imbalance markets, see, for example, \cite{jiang2015optimal} or \cite{lohndorf2023value}.}
\end{itemize}

We consider an optimal trading problem with a general propagator model under linear functional inequality constraints that subsumes and extends the two motivating applications in a stochastic setting and continuous time. \\

\begin{figure}[H]
    \begin{center}
    \includegraphics[scale=0.45]{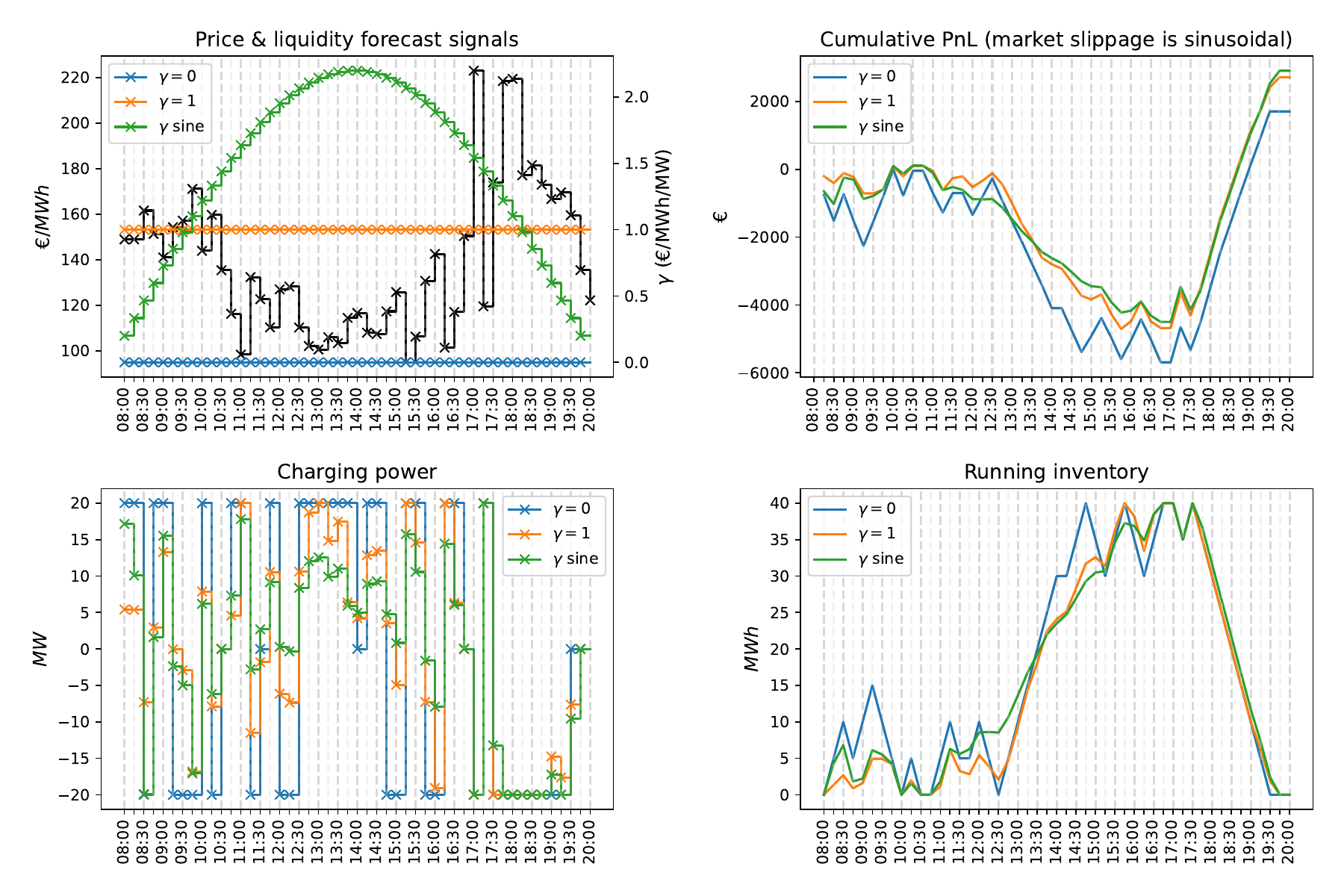}
    \end{center}
    \vspace{-0.4cm}
    \caption{Motivating example: battery optimization in deterministic discrete time \eqref{eq:quadratic_program_functional_with_transient}--\eqref{eq:first_example_linear_program_constraints} solved as a Linear Program in the case $\gamma_i = 0, \; i \in \{0,\ldots, N-1 \}$ (in blue) and as Quadratic Programs for two specifications of the slippage costs $\gamma_i > 0, \; i \in \{0,\ldots, N-1 \}$ (constant in orange and sinusoidal in green in the plots) and zero weights $K_{i,j} = 0, \; i,j \in \{0,\ldots, N-1 \}$. Top-left: discrete time price (in black) and slippage costs forecasts (assumed to be piecewise constant on each quarter); top-right: running cumulative PnL \eqref{eq:quadratic_program_functional_with_transient} assuming the market slippage cost intensities are sinusoidal (green curve in top-left plot); bottom-left: optimal charging powers (with a standard deviation of $18.9$ without slippage costs, $13.5$ for the constant slippage costs and $12.4$ for the sinusoidal one); bottom-right: optimal running load capacities. Constraints are taken as: $(u^{min}, u^{max}) = (-20, 20)$ (in $MW$) and $(X^{min}, X^{max}) = (0, 40)$ (in $MWh$) in \eqref{eq:first_example_linear_program_constraints}.}\label{F:Quadratic_program_example}
\end{figure}

To motivate our class of problems, we describe a deterministic version in discrete time. Assume that we have access to discrete time price forecasts \(S := \left( S_{i} \right)_{i \in \{0, \cdots, N\}}\) on a uniform time grid $0=t_0<t_1<\ldots<t_N$, where \(N \in \mathbb{N}^{*}\). These could be quarter hour electricity prices (illustrated in black in the upper-left plot of Figure~\ref{F:Quadratic_program_example}). Our goal is to determine the optimal trading (or charging)  strategy $u := (u_{i})_{i\in \{0,\ldots, N-1 \}}$  while adhering to constraints on both $u$ and on the inventory $X^u$ (or capacity of battery), defined by 
$$   X_{i}^{u} := X_{0} + \Delta \sum_{j=0}^{i-1} u_{j}, \quad i \in \{0, \cdots, N\},$$
with $X_0$ the initial inventory and $\Delta$ the time step. Inspired by the discrete time propagator model of price impact from \cite{bouchaud2009price}, we aim to maximize the final cumulative Profit and Loss (PnL) over the trading horizon:
\begin{align}\label{eq:quadratic_program_functional_with_transient}
    J(u) & := \underbrace{- \sum_{i=0}^{N-1} \Big( S_{i} + \frac{\gamma_{i}}{2} u_{i} + \sum_{j=0}^{i-1} K_{i,j} u_{j} \Big) \Delta u_{i}}_{\text{Cumulative cash-flows with market impact}} + \underbrace{X_{N}^{u} S_{N}}_{\text{Final inventory valuation}},
\end{align}
subject to the constraints
\begin{equation}\label{eq:first_example_linear_program_constraints}
    \begin{cases}
        u_{i}^{\min} \leq u_{i} \leq u_{i}^{\max} , \quad  X_{i+1}^{\min} \leq X_{i+1}^{u} \leq X_{i+1}^{\max},  \quad i \in \{0, \cdots, N-1\}, \\
        X_{0} = X_{N}^{u} = 0.
    \end{cases}
\end{equation}
The deterministic positive terms $(\gamma_i)_{i \in \{0,\ldots, N-1 \}}$ represent \textit{the slippage costs i.e. the instantaneous market impact} (or \textit{the running operating costs} of the asset) while the non-negative weights $\left( K_{i,j} \right)_{i,j \in \{0,\ldots, N-1 \}}$ take into account \textit{the transient market impact of past trades} on the price. Since the optimization problem is discrete and deterministic, it can be reformulated as a Linear Program (LP) when $\gamma_i = 0, \; i \in \{0,\ldots, N-1 \}$ and as a Quadratic Program (QP) as soon as $\gamma_i > 0, \; i \in \{0,\ldots, N-1 \}$ in the sense of \cite[Chapter~4, Sections~4.3 and 4.4]{boyd2004convex}, and solved numerically using commercial or open-source solvers\footnote{\label{note1} We used CVXOPT in our case, see \hyperlink{https://cvxopt.org}{https://cvxopt.org}, which employs interior-point methods, see \cite{vandenberghe2010cvxopt} and references therein.}. Figure~\ref{F:Quadratic_program_example} displays the resulting optimal strategies as well as the associated running inventories for different specifications of slippage costs and without transient impact, i.e.~\( \left( K_{i,j} \right)_{0 \leq j < i \leq N-1} = 0 \), and the respective cumulative cash flows assuming the market slippage costs are sinusoidal given by the green curve in the top-left plot.

\vspace{-0.2cm}

When comparing strategies with zero (in blue), positive constant (in orange) and sinusoidal (in green) slippage cost curves, notice that the lower the slippage costs, the more aggressively the agent trades, with the limit case of a `bang-bang' strategy when there are no slippage costs (in blue). Furthermore, when assuming the realized market slippage costs are indeed sinusoidal, notice that accounting for a constant slippage cost (orange curves) yields a higher cumulative PnL at maturity than neglecting them (blue curves). 

A major limitation of the above model is the deterministic nature of the signal $(S_i)_{i \in \{0,\ldots, N-1 \}}$. A stochastic signal is crucial for modeling real-world uncertainties, to make the problem formulation more realistic and applicable. However, randomness complicates both theoretical analysis and numerical solutions, especially with non-Markovian processes. Even in discrete time, the optimization  problem can no longer be reformulated as a standard LP or QP problem.

\paragraph{Contributions.} The main contributions of this paper is to extend and solve the control problem formulation \eqref{eq:quadratic_program_functional_with_transient}--\eqref{eq:first_example_linear_program_constraints}
\begin{enumerate}
    \item[(i)] in the continuous-time setting (i.e.~when $\Delta \to 0$),
    \item[(ii)] where the forecast signal price and the constraining functions are general stochastic processes, and not necessarily Markovian,
    \item[(iii)] with both equality and inequality constraints linear in the trading rate and the associated inventory,
    \item[(iv)] with a general linear propagator model for modelling the  transient part of the market impact,
    \item[(v)] by generalizing the well-known Uzawa algorithm in order to build admissible trajectories of the stochastic optimal control.
\end{enumerate}

 From the theoretical side, we prove the well-posedness of the problem, formulated as a stochastic convex programming problem, even when dealing with a general exogenous price and stochastic constraints. We derive necessary and sufficient infinite-dimensional optimality conditions of Karush-Kuhn-Tucker (KKT) type \cite{karush1939minima, kuhn2013nonlinear}, following the approach of \cite{rockafellar1976stochastic_3} and leveraging the generalized Lagrange multiplier rule. We explicitly solve the first-order conditions expressed as a stochastic Fredholm equation of the second kind. The optimal control is explicitly given in terms of the corresponding Lagrange multipliers associated with the constraints and their conditional expectations.
 Our main results, summarized in  Theorems~\ref{T:generalized_kkt_theorem_necessary} and \ref{T:generalized_kkt_theorem_sufficient}, are valid in a wide class of real-world problems beyond the typical Markovian assumptions.

From the numerical perspective, we propose in Section~\ref{S:numerics} a stochastic formulation of the Uzawa algorithm on the dual problem to construct the stochastic Lagrange multipliers numerically via a stochastic projected gradient ascent, combined with a least-squares Monte Carlo regression step to approximate their conditional expectations.
We illustrate our findings, in Section~\ref{S:Applications}, on two practical applications with stochastic signals: (i) an optimal execution problem with transient impact, with either a `no-shorting' constraint in the presence of a `sell' signal or a `no-buying' constraint in the presence of a `buy' signal and a stochastic `stop-trading' constraint whenever the exogenous price drops below a specified reference level; (ii) a battery storage problem with constant running operating costs, seasonal signals, and fixed constraints on both the charging power and the load capacity of the battery.

\paragraph{Related literature.} 
The stochastic target problem in optimal execution is studied in \cite{ackermann2024reducing}, where the final inventory level is a prescribed random variable, with quadratic running cost from a given running target level while in \cite{aksu2023optimal}, an Almgren-Chriss control problem is considered where the full liquidation constraint is relaxed if the price goes below a specified threshold. The impact of `no-shorting' constraints has been explored on portfolio construction in \cite{li2002dynamic, cui2014optimal} as well as on asset prices and liquidity provision in \cite{bai2006asset}, while a participation constraint is tackled in \cite{labadie2012optimal} in the context of target close and implementation shortfall trading algorithms. Also, `no direction change' constraints on the trading rate have been explored in the context of $N$-player and mean-field games for optimal portfolio liquidation in \cite{fu2024mean}, which are similar to the `no buy' constraint tackled in this paper.

Several works have studied the modeling of energy storage systems such as batteries or dams either by simplifying the problem in discrete time or by imposing some Markovian restrictions on the battery's load capacity such that the asset must discharge when it is full, and charge when it is empty, see for example \cite{collet2018optimal}, \cite{ben2023lagrangian}, \cite{leclere2014contributions} or using Deep Learning techniques \cite{lemaire2024swing}. 

The modelling of continuous-time inequality constraints has already been explored by \cite{chiarolla2013generalized} where they solve the $N$-firm stochastic irreversible investment problem under limited resources: they propose a stochastic extension of the KKT theorem and even provide some explicit solutions for the Lagrange measure associated to their resource constraint. More recently, \cite{bayraktar2019controlled} solve a stochastic control problem with path-dependent coefficients under deterministic inequality constraints on a state variable dependent on the control, by solving a backward stochastic partial differential equation with Neumann boundary conditions satisfied by the value function.

From the theoretical point of view, our work greatly echoes with all the contributions brought by Rockafellar and Wets in the 1970's on a two-stage stochastic convex program with recourse, where the dual problem formulation is extensively explored via perturbation theory, in particular in their works \cite{rockafellar1976stochastic_1, rockafellar1976stochastic_3, rockafellar1976stochastic_2, rockafellar1975stochastic_0, rockafellar1977measures_4}. The Economic theory literature has also demonstrated a great interest in the dual methods in Economy equilibrium and growth models, where sufficient conditions are derived to be able to interpret the Lagrange multipliers as \textit{price processes}, i.e.~such that their singular part in their Yosida-Hewitt decomposition \cite{yosida1952finitely} is zero, see for example \cite{bewley1972existence}, \cite{morand2015lagrange}.

Originally formulated in the deterministic finite-dimensional setting in \cite[p.154-165]{arrow1958studies}, the \textit{Uzawa algorithm} has been extensively used and studied in the literature of constrained optimization and saddle-points approximations, and has been extended in the infinite-dimensional setting on Hilbert spaces, see for example the bibliographical comments in \cite[Section 1.2]{bertucci2020remark}, and even in Banach spaces \cite{leclere2014contributions}. Extending the Uzawa algorithm into the stochastic setting is not new with, for example, some convergence results already obtained in \cite{leclere2014contributions} for the case of equality constraints and some initial numerical results presented in \cite{barty2010decomposition} in the Markovian case, where the primal problem is solved using dynamic programming. 

Finally, our modelling framework for both the temporary and transient market impacts is greatly inspired by the works of \cite{gatheral2012transient}, \cite{neuman2022optimal}, \cite{abi2022optimal} and in particular the work \cite{jaber2023equilibrium} where stochastic Fredholm equations of the second kind have been introduced and solved explicitly for the first time.

\paragraph{Outline.} The paper is organized as follows. In Section~\ref{S:problem_formulation}, we formulate the general stochastic control problem of propagator type under almost sure constraints and we present our main theoretical results. In Section~\ref{S:numerics}, we define a stochastic extension of the Uzawa algorithm to obtain numerically admissible sample trajectories of the optimal stochastic control. In Section \ref{S:Applications}, we apply the algorithm to an optimal execution problem and  an optimal storage management under almost sure constraints. Finally, we detail in Sections~\ref{S:existence} and \ref{ss:proof_theorem_necessary_sufficient_condition} the proofs of our main results. 

\paragraph{Notations.}
We fix a finite time horizon $T > 0$ and a filtered probability space $(\Omega, \mathcal F, (\mathcal F_t)_{t \in [0,T]}, \mathbb{P})$ satisfying the usual conditions, where $\Omega$ is Polish\footnote{One can think of $\Omega$ as the Skorohod space for example, see \cite[Chapter 3]{billingsley2013convergence}.} and $\mathbb{P}$ is a probability measure. We denote the Lebesgue measure by $\d t$ on the Borel sigma-algebra $\mathcal{B} \left( [0,T] \right)$, and we denote by $\d t \otimes \mathbb{P}$ the product measure on the sigma-algebra $\mathcal{B} \left( [0,T] \right) \otimes \mathcal{F}$.
We introduce the standard spaces
\begin{equation} \label{eq:def_mathcal_L_p}
    \mathcal{L}^{p} := \left\{ f: [0,T] \times \Omega \to \R \; \text{prog.~measurable}, \; \int_{0}^{T}\E \left[  |f_{t}|^{p}\right] \d t  < \infty \right\},  \quad p \geq 1. 
\end{equation}
For $p=2$,  we equip $\mathcal L^2$ with the inner-product
\begin{equation}
    \langle f, g \rangle_{\mathcal{L}^{2}} := \int_{0}^{T}\E \left[  f_{t} g_{t}  \right]\d t, \quad f,g \in \mathcal{L}^{2},
\end{equation}
which makes it a Hilbert space with the  associated norm $\| f\|_{_2} := \sqrt{\langle f,f\rangle_{\mathcal{L}^{2}}}$. Similarly, $L^{2}$ denotes the standard Hilbert space of square-integrable $[0,T]$-valued functions with respect to the Lebesgue measure, with its inner product denoted by $\langle ., . \rangle_{L^{2}}$. Also, for $p=\infty$,  we consider 
\begin{equation} \label{eq:def_mathcal_L_infty}
    \mathcal{L}^{\infty} := \left\{ f: [0,T] \times \Omega \to \R \; \text{prog.~measurable}, \; \exists C \in \mathbb{R}_{+}, \; \left| f_{t}(\omega) \right| \leq C \text{ hold } \left( \d t \otimes d\mathbb{P} \right)-a.e. \text{ on } [0,T] \times \Omega \right\},
\end{equation}
which is a Banach space with the norm
\begin{equation}
    \left| \left| f \right| \right|_{\infty} := \inf_{C \in \mathbb{R}_{+}} \left\{ \left| f \right| \leq C \text{ hold } \left(\d t \otimes \P \right)-a.e. \text{ on } [0,T] \times \Omega \right\}, \quad f \in \mathcal{L}^{\infty}.
\end{equation}
For the sake of conciseness, we abusively use the notation $\langle \cdot,\cdot  \rangle$ either to refer to the dual product between $\left( \mathcal{L}^{\infty} \right)^{p}, \; p \geq 1$ and its dual space, or to the standard Euclidean scalar product in $\mathbb{R}^{k}$, for some $k \in \mathbb{N}$; there should be no ambiguity based on the context. Finally, for $t \in [0,T]$, we denote by $\E_t$ the conditional expectation $\E[\cdot | \mathcal F_t]$.

\section{Problem formulation and main theoretical results} \label{S:problem_formulation}

\subsection{Problem formulation}
 
\noindent We consider an agent with an initial holding $X_0 \in \R$ controlling the rate of change $u$ of her inventory such that
\begin{equation} \label{eq:inventory_def}
    X_{t}^u := X_{0} + \int_0^t u_s \d s, \quad t \in [0,T].
\end{equation}
The processes $X^{u}$ and $u$  can represent: 
\begin{enumerate}
    \item[(i)] the inventory $X^u$ held in a certain asset together with the trading rate $u$ (in $s^{-1}$) in an execution problem, such that $u>0$ means the agent buys the asset, while $u<0$ means the agent shorts it;

    \item[(ii)] the  capacity of a battery $X^u$ (typically in $MWh$) and the charging power $u$ (typically in $MW$) in a battery storage problem.
\end{enumerate}

\noindent Given an exogenous price process $S$ of the form 
\begin{equation} \label{eq:exogenous_price_process}
    S_t=P_{t}+M_t,
\end{equation}
with $M$ a centered square-integrable martingale and   $P$  a finite variation process in $\mathcal L^2$,  the agent aims at minimizing the following functional associated to a control $u$ over a finite time horizon $T$:
\begin{equation} \label{eq:functional_formulation_storage}
    \mathcal{J}(u) := \mathbb E \left[ \int_{0}^T \left( S_t + \frac{\gamma(t)}{2} u_t + \int_0^t K(t,s)u_s \d s \right) u_t \d t - X^u_T S_T \right],
\end{equation}
where $\gamma$ is a measurable and bounded positive function with $\inf_{[0,T]} \gamma > 0$, and  $K$ $:$ $[0,T]^2 \to \R_{+}$ is an admissible  kernel  in the sense of Definition~\ref{Def:admissiblekernel} below which will ensure that the minimization problem is well-posed. 

\begin{remark}[Interpretation of the functional \eqref{eq:functional_formulation_storage}] The \textit{effective} exchange price process defined by 
    \begin{equation}
        \tilde{S}_t := S_t + \frac{\gamma(t)}{2} u_t + \int_0^t K(t,s) u_s \d s, \quad t \in [0,T],
    \end{equation}
    is the price at which the agent trades the asset while facing both linear market slippage (or running operating) costs with intensity $\gamma$ as well as \textit{transient price impact}
    \begin{equation} \label{eq:transient_market_impact_state_variable}
        Z_{t}^{u} := \int_{0}^{t} K(t,s) u_{s} \d s, \quad t\leq T,
    \end{equation}
    due to her past interaction with the market weighted by the kernel $K$. Hence, the term $ \int_0^T \tilde S_t \d X^u_t = \int_0^T \tilde S_t u_t \d t  $ that appears in $\mathcal J(u)$ corresponds to the profit and loss of trading.  The last term $X^u_T S_T$ in the functional stands for the valuation of the final inventory  $X^u_T$  multiplied by the exogenous price $S_T$ at terminal date $T$; note that we do not use the effective price for the final valuation, in order to avoid fake profits optimization coming from the impact model itself, see \cite{kolm2023you}. 
    Note that the price process $S$ from \eqref{eq:exogenous_price_process} not only represents the exogenous price but can also incorporate a stochastic trading signal (which can be a forecast of future prices) in the drift process $P$. This model, known as the propagator model,  was initially introduced in \citet{bouchaud2009price}, \citet[Chapter 13]{bouchaud2018trades}, and has been extensively studied since then, see \citet*{abi2022optimal, abijaberOptPort2024,benzaquen2017dissecting,gatheral2012transient,graewe2017optimal,neuman2022optimal,obizhaeva2013optimal}.
\end{remark}

\begin{definition}\label{Def:admissiblekernel}
    We define the set $\mathcal K$ of admissible kernels by:  $K\in \mathcal K$ if 
\begin{itemize}
    \item[(i)] $K$ is a measurable square-integrable \textit{Volterra} kernel, i.e.
    \begin{equation} \label{eq:volterra_kernel}
        K(t,s) = 0, \quad 0 \leq t < s \leq T,
    \end{equation}
    and 
       \begin{equation} \label{eq:square_integrable_kernel_assumption}
        \int_{0}^{T} \int_{0}^{T} K(t,s)^{2} \d s \d t < \infty.
    \end{equation}
    \item[(ii)] $K$ is \textit{non-negative definite}, i.e.~for all $f\in L^2$,
    \begin{equation} \label{eq:pos_def}
        {\int_{0}^T \int_{0}^T K(t,s) f(s)f(t)\d s\d t \;  \geq \;  0.
        }
    \end{equation}
    \end{itemize}
\end{definition}

\begin{example} \label{E:kernels}
Two notable admissible kernels are given by the exponential kernel\footnote{Any linear combination of exponential kernels is also admissible, to be able to deal with multi-scale exponential decays.}
\begin{equation} \label{eq:def_exponential_kernel}
    K(t,s) := ce^{-\rho (t-s)}\mathbb{1}_{\{t > s\}}, \quad t,s \in [0,T], \quad c, \; \rho > 0,
\end{equation}
and the fractional kernel 
\begin{equation} \label{eq:def_fractional_kernel} 
    K(t,s) := c(t-s)^{\alpha-1} \mathbb{1}_{\{t > s\}}, \quad t,s \in [0,T], \quad c>0, \quad \alpha \in \big( 0,1 \big).
\end{equation}
 We refer to \cite[Example~2.5]{abi2022optimal} for additional examples.
\end{example}

The final ingredient for our problem formulation is the introduction of four constraining processes $(a^{i})_{i \in \{1, 2, 3, 4 \}} \in \left(\mathcal{L}^{\infty}\right)^{4}$ such that 
\begin{equation} \label{eq:def_constraining_processes}
    a^{1} \leq a^{2}, \quad a^{3} \leq a^{4} \; \text{ and } \; a_{0}^{3} \leq X_{0} \leq a_{0}^{4}, 
\end{equation}
and the  set of admissible controls  defined by
\begin{equation} \label{eq:def_admissible_set_of_controls}
    \mathcal{U} := \left\{ u \in \mathcal{L}^{\infty}, \; \; a_{t}^{1} \leq u_{t} \leq a_{t}^{2} \; \text{and} \; a_{t}^{3} \leq X_{t}^{u} \leq a_{t}^{4} \text{ hold } \left( \d t \otimes \mathbb{P} \right) -a.e. \text{ on } [0,T] \times \Omega \right\}.
\end{equation}

\noindent The agent aims to solve the following constrained optimization problem
\begin{equation} \label{eq:problem_formulation}
    \inf_{u \in {\cal U}} \mathcal{{J}}(u).   
\end{equation}

Our framework encompasses several practical configurations for optimal trading in quantitative finance and energy storage management as shown by  the next examples.

\begin{example}[Optimal trading with  constraints] \label{ex:stochastic_target}
    In the case of optimal trading, the set $\mathcal U$ is versatile enough to  represent several  constraints of regulatory or managerial type for the trading rate $u$ and the associated inventory $X^{u}$ such as:
    \begin{enumerate}
        \item[(i)] deterministic or stochastic flow: $X_T^u=b_T$ for some random variable $b_T$, the case $b_T=0$ corresponds to the total liquidation constraint, whenever $a^3_T=a^4_T=b_T$,
        \item [(ii)] no buying constraint: $u_t\leq 0$, for all $t\in [0,T]$, whenever $a^2 = 0$,
        \item [(iii)] no shorting constraint:  $X^u_t\geq 0$, for all $t\in [0,T]$ whenever  $a^3 = 0$,
        \item [(iv)] index/portfolio tracking and stochastic target problems: $P_t - \epsilon \leq X_t\leq P_t + \epsilon$, where $P$ is a certain stochastic reference portfolio and $\epsilon >0$, whenever $a^3 = P - \epsilon$ and $a^4 = P + \epsilon$.
        \item[(v)] Conditions on minimum price in the spirit of \cite[Section 2.1]{aksu2023optimal}: 
       \begin{itemize}
           \item[a)] the total execution of $X_{0}>0$ shares happens only if the closing fundamental price \eqref{eq:exogenous_price_process} does not drop below a certain reference price $S^{r}>0$, else relax the full liquidation constraint to an inventory level of $\delta > 0$ shares at maturity, whenever $a^3 = 0$ and $a^4_T = \delta \mathbb{1}_{(-\infty, S^{r}]}(S_{T})$,
           
           \item[b)] liquidating the $X_{0}>0$ shares is allowed until the price \eqref{eq:exogenous_price_process} goes below a given reference price $S_{0}>S^{r}>0$ at which point trading is completely stopped, and the position is no longer constrained to be fully closed: whenever $- a_{.}^{1} = a_{.}^{2} = \bar M \mathbb{1}_{[0, \tau^{r}]}(.)$, $- a_{t}^{3} = a_{t}^{4} = \bar M, \; t \in [0,T)$, and $- a_{T}^{3} = a_{T}^{4} = \bar M \mathbb{1}_{[0, T)}(\tau^{r})$, for $\tau^{r} := \inf_{t \in [0,T]} \left\{ S_{t} < S^{r} \right\} \wedge T$ and for some prescribed large enough constant\footnote{As argued more in detail in Section \ref{ss:optimal_execution_under_constraints}, $\bar M$ is taken large enough such that everything happens as if there is no constraint on the resulting optimal numerical control's sample trajectories.} $\bar M>0$, whose unit is obvious depending on the constraint.
       \end{itemize}
       Examples combining $(i)$ with either $(ii), (iii)$ or $(v)-b)$ will be solved numerically in Section \ref{ss:optimal_execution_under_constraints} in the context of optimal execution.
    \end{enumerate}
\end{example}

\begin{example}[Energy storage management] \label{ex:battery_swing_execution}
    In the case of optimal battery management, the set $\mathcal U$ conveniently represents physical constraints for the charging power $u \in \mathcal{L}^{\infty}$ and the associated load capacity $X^{u}$ that typically read as 
    \begin{align} \label{eq:charging_constraints}
        -u^{max} \leq u_t \leq u^{max}, \quad   0 \leq X^{u}_t \leq X^{max}, \quad t \in [0,T], 
    \end{align}
    where the constant bounds $u^{max}, X^{max} > 0$ are respectively interpreted as the maximum power delivered by the battery and its maximum load capacity.
    Notice that the control $u \equiv 0$ always lies in $\mathcal{U}$, which is consequently non-empty. Then, solving problem \eqref{eq:problem_formulation}  amounts to optimizing the expected realized PnL of a battery asset connected to a real-time electricity market. A specification of such optimal storage problem with running operating costs and seasonal stochastic signals is solved numerically in Section \ref{ss:energy_storage_example}.
\end{example}

\begin{remark}[Why choosing $\mathcal{L}^{\infty}$?] \label{R:choice_L_infty}
    First and foremost, from a practical point of view, it seems reasonable to assume that, whatever the state of the world $\omega \in \Omega$ at observation date $t \in [0,T]$, the agent's trading rate (or charging power) $u_{t}(\omega)$ will always be bounded due to technological restrictions\footnote{What does an infinite trading rate mean anyway?}.
    
    From the theoretical point of view, the choice of $\mathcal{L}^{\infty}$ space for the admissible controls in \eqref{eq:def_admissible_set_of_controls} is convenient to solve the constrained control problem \eqref{eq:problem_formulation} since its associated positive cone, defined by \begin{equation} \label{eq:positive_cone_L_infty}
    \mathcal{L}_{+}^{\infty} := \left\{ f \in \mathcal{L}^{\infty}, \; f \geq 0 \text{ hold } \left( \d t \otimes \d\P \right) - a.e. \right\},
    \end{equation} 
    has a non-empty interior, contrary to all respective positive cones associated to the other standard $\mathcal{L}^{p}$ Banach spaces for $1 \leq p < \infty$. Such property is of paramount importance to obtain the existence of Lagrange multipliers via a geometric version of the Han-Banach separation Theorem recalled in Section~\ref{ss:existence_lagrange_linear_form}, see, for example, \cite[Section 8.2, p.219]{luenberger1997optimization}, \cite[Section 4.1.2]{leclere2014contributions} or the contributions by \citet{rockafellar1975stochastic_0, rockafellar1976stochastic_1, rockafellar1976stochastic_2, rockafellar1976stochastic_3} for additional background on the standard choice of the essentially bounded functions in constrained optimization. Although $\mathcal{L}^{\infty}$ does not enjoy neither separability nor reflexivity, we will leverage the fact that
    \begin{equation}
        \mathcal{L}^{\infty} = \left( \mathcal{L}^{1} \right)^{\star},
    \end{equation}
    in Subsection~\ref{ss:existence_uniqueness} to prove the existence of a solution to \eqref{eq:problem_formulation}, see, for example, \cite[Theorem 4.14]{brezis2011functional}.
\end{remark}

\subsection{The solution of  the control problem} \label{SS:main_results}

 Our main Theorems~\ref{T:generalized_kkt_theorem_necessary} and \ref{T:generalized_kkt_theorem_sufficient} below provide the unique solution $\hat u$ to the constrained optimization problem \eqref{eq:problem_formulation} explicitly in terms of stochastic Lagrange multipliers. Before stating them, we need to introduce some additional notations and recall some preliminary results from functional analysis.\\

The dual space $(\mathcal{L}^{\infty})^{\star}$ of $\mathcal{L}^{\infty}$ is of special interest in this paper, and is identified with the space $\textbf{\textit{ba}}$ of bounded finitely-additive set-valued functions on the product sigma algebra $\mathcal{B}\left( \left[ 0,T \right] \right) \otimes \mathcal{F}$ which vanish on sets of $\left( \d t \otimes \P \right)$-measure zero. 
The set $\textbf{\textit{ba}}$ is sometimes also referred to as \textit{finitely additive measures} following the terminology\footnote{Such terminology can cause some confusion since by today's standard definition, a measure is a set-valued function that is countably additive, which is not necessarily the case for elements of $\textbf{\textit{ba}}$, as they are only required to be finitely additive.} of \citet{yosida1952finitely}. For any element $\psi: \mathcal{B}([0,T]) \otimes \mathcal{F} \to \R$ of the $\textbf{\textit{ba}}$ space, {we consider the dual pairing
\begin{equation} \
    \langle u, \psi \rangle := \int_{[0,T] \times \Omega} u_{t}(\omega) \psi(\d t,\d \omega), \quad  u \in \mathcal{L}^{\infty},
\end{equation}
where the right-hand side is a \textit{Radon's integral}, see for example \cite[Chapter IV, Section 8, Theorem 16]{dunford1988linear} for a detailed construction of such an integral, and $\textbf{\textit{ba}}$ is normed by
\begin{equation}
    ||\psi||_{\textbf{\textit{ba}}} \; := \sup_{||u||_{\infty} \leq 1} \langle u, \psi \rangle, \quad \psi \in \textbf{\textit{ba}}.
\end{equation}}
We also denote by $\textbf{\textit{ba}}_{+}$ the set of non-negative elements in $\textbf{\textit{ba}}$ defined by 
\begin{equation} \label{eq:non_negativity_condition}
  \textbf{\textit{ba}}_{+} := \left\{ \psi \in \textbf{\textit{ba}}: \; \langle u, \psi \rangle \geq 0, \;   u \in \mathcal{L}_{+}^{\infty} \right\}.
\end{equation}
Additionally, we recall the Yosida-Hewitt decomposition which states that any finitely additive measure can be uniquely decomposed as the sum of \textit{a countably additive part} and \textit{a purely finitely additive part}, such that for any $\psi \in \textbf{\textit{ba}}$
\begin{equation} \label{eq:hewitt_yosida_decomposition}
    \langle u, \psi \rangle = \int_{[0,T] \times \Omega} u_{t}(\omega) \lambda_{t}(\omega) \d t \P(\d \omega) + \psi^{o}(u), \quad  u \in \mathcal{L}^{\infty},
\end{equation}
where $\lambda \in \mathcal{L}^{1}$ is the density of the countably additive part, with $\mathcal{L}^{1}$ defined by \eqref{eq:def_mathcal_L_p}, and $\psi^{o} \in \mathcal{S}$, $\mathcal{S}$ being the set of continuous linear functionals $\psi^{o}$ on $\mathcal L^{\infty}$ such that there exists an increasing sequence of measurable sets $A_{k}$ such that $\bigcup_{k \geq 0} A_{k} = [0,T] \times \Omega$, and for each $k$ one has $\psi^{o}(v)=0$ for all the functions $v \in \mathcal{L}^{\infty}$ vanishing almost everywhere outside of $A_{k}$, see \cite[Theorem 1.24]{yosida1952finitely} for the original decomposition formulation or \cite[Equation (2.4)]{rockafellar1976stochastic_3} where such decomposition is used in the context of stochastic programming. If additionally $\psi \in \textbf{\textit{ba}}_{+}$, then $\lambda \geq 0$ and $\psi^{0} \in \mathcal{S}_{+} := \mathcal{S} \cap \textbf{\textit{ba}}_{+}$ by \cite[Theorem 1.23]{yosida1952finitely}.\\

 In a nutshell, our methodology for solving the optimization problem  \eqref{eq:problem_formulation} in Theorem~\ref{T:generalized_kkt_theorem_necessary} consists in two main steps that we briefly sketch now without much rigour.\\

\textbf{First step.} We recast our optimization problem \eqref{eq:problem_formulation} as an infinite-dimensional convex problem on the space $\mathcal L^{\infty}$ by writing 
\begin{align} \label{eq:stochasticConvexProgramFormulation}
     \mathcal{J}(u) & = - \E \left[ \int_{0}^T \left(\alpha_t - \left( \frac{\gamma(t)}{2} u_{t} +\int_0^t K(t,s)u_{s} \d s \right) \right) u_{t} \d t \right] - X_{0} \E \left[ S_T \right], \\
    \mathcal{U} & = \left\{ u \in \mathcal{L}^{\infty}, \; \mathcal{G}(u) \in -\mathcal{C} \right\},
\end{align}
where 
\begin{equation} \label{eq:alpha_signal_definition}
     \alpha_{t} := \E_{t} \left[ P_{T} - P_{t} \right], \quad t \in [0,T],
\end{equation}
$\mathcal{G}$ is a linear constraining functional such that
\begin{equation} \label{eq:def_constraining_functional}
    \mathcal{G}(u):=\left(a^{1}-u, u-a^{2},  a^{3} - X^u, X^u - a^{4} \right), \quad u \in \mathcal{L}^{\infty},
\end{equation}
and $\mathcal{C} := \left(\mathcal{L}_{+}^{\infty}\right)^{4}$ is a non-negative ordering cone with non-empty interior. We  define the associated Lagrange functional
\begin{equation} \label{eq:def_lagrange_functional}
    \mathcal{L} \left( u, \psi \right) := \mathcal{J}(u) + \langle \mathcal{G}(u), \psi \rangle, \quad u \in \mathcal{L}^{\infty}, \quad \psi \in \textbf{\textit{ba}}_{+}^{4}. 
\end{equation}

\noindent The above convex program admits a unique solution $\hat{u}$, and by applying the generalized Lagrange multiplier rule (see for example \citet[Chapter 5]{jahn1994introduction}, or \citet[Chapter 8]{luenberger1997optimization}) we obtain the existence of four non-negative, bounded and finitely additive set-valued functions $\psi :=(\psi^{1},\psi^2,\psi^3,\psi^4) \in \textbf{\textit{ba}}_{+}^4$ satisfying
\begin{itemize}
\item [\textbf{(i)}] a Lagrange stationary equation in the form 
\begin{align}\label{eq:lagstation}
    \nabla_u {\mathcal{J}}(\hat u ) + \langle \nabla_u \mathcal{G}(\hat u ), \psi \rangle= 0 
\end{align}
where $\nabla_u$ denotes the Fréchet derivative, 
    \item [\textbf{(ii)}]
     and such that the respective densities of the finitely additive parts of $\left(\psi^{i}\right)_{i \in \{1,2,3,4\}}$ from \eqref{eq:hewitt_yosida_decomposition}, denoted by $\left(\lambda^{i}\right)_{i \in \{1,2,3,4\}} \in \left( \mathcal{L}_{+}^{1} \right)^{4}$, satisfy the complementary slackness conditions given by
\begin{equation} \label{eq:KKT_slackness_conditions}
    \begin{cases} 
       0 = \mathbb{E} \left[ \int_{0}^{T} \left(a_{t}^{1} - \hat{u}_{t}\right) \lambda_{t}^{1} \d t \right]  \\
       0 = \mathbb{E} \left[ \int_{0}^{T} \left(\hat{u}_{t} - a_{t}^{2}\right) \lambda_{t}^{2} \d t \right]  \\
       0 = \mathbb{E} \left[ \int_{0}^{T} \left( a^3_t  - X^{\hat{u}}_t \right) \lambda_{t}^{3} \d t \right] \\
       0 = \mathbb{E} \left[ \int_{0}^{T} \left( X^{\hat{u}}_t - a^4_t \right) \lambda_{t}^{4} \d t \right]
    \end{cases}.
\end{equation}
\end{itemize}
\noindent The four stochastic processes $\left(\lambda^{i}\right)_{i \in \{1,2,3,4\}}$ play the role of infinite-dimensional Lagrange multipliers associated with the four respective inequality constraints introduced in $\mathcal{U}$. Interestingly, the singular parts of $\psi$, recall \eqref{eq:hewitt_yosida_decomposition},  can be discarded as they do not enter in the expression of the optimal control $\hat u$ as shown in the sequel.\\

\textbf{Second step.} We explicitly solve the Lagrange stationary equation \eqref{eq:lagstation} to obtain the optimal control $\hat{u}$. In the following, we set $\gamma = 1$ in \eqref{eq:functional_formulation_storage} without loss of generality, i.e.~it suffices to replace $K$ by $K/\gamma$ and $S$ by $S/\gamma$. Then, Equation \eqref{eq:lagstation} leads to a stochastic Fredholm equation of the second kind of the form
\begin{equation} \label{eq:Fredholm_equation}
    \hat{u}_{t} + \int_0^{t} K(t,s) \hat{u}_{s} \d s + \int_t^T K(s,t) \E_{t} \hat{u}_{s} \d s = \varphi_{t} + \alpha_{t} \quad \left(\d t \otimes \P\right) - a.e.,
\end{equation}
where $\alpha$ is given by \eqref{eq:alpha_signal_definition}, and $\varphi \in \mathcal{L}^{2}$ satisfies the important decomposition\footnote{Given the functional \eqref{eq:stochasticConvexProgramFormulation} to minimize is an expectation, the existence of an aggregated countably additive Lagrange multiplier absolutely continuous with respect to $dt \otimes \mathbb{P}$ with density $\varphi$ given by \eqref{eq:aggregated_lagrange_multiplier} and associated to the constraining functional $\mathcal{G}$ from \eqref{eq:def_constraining_functional} is not surprising, see for example \cite{bewley1972existence} and \cite{morand2015lagrange}[Section 3] where they study sufficient conditions to cancel the singular parts of the multipliers for general stochastic convex programs.}
\begin{equation} \label{eq:aggregated_lagrange_multiplier}
     \varphi_{t} = \lambda^{1}_t - \lambda^{2}_t + \mathbb{E}_{t} \left[ \int_{t}^{T} \left( \lambda_{s}^{3} - \lambda_{s}^{4} \right) \d s \right].
\end{equation}
Equation~\eqref{eq:Fredholm_equation} was recently  solved in \citet*[Section 5]{jaber2023equilibrium} and further studied in \citet*{abijaberOptPort2024}.  We first introduce some notations. For any admissible kernel $K \in \mathcal K$, we will use the following notation
\begin{equation}
    K_{t}(s,r) := K(s,r) \mathbb{1}_{r \geq t}, \quad t \in [0,T],
\end{equation}
and we will use the bold notation $\bm{K} : L^{2} \to L^{2}$ to refer to the bounded linear operator induced by a square-integrable kernel $K$ defined by
\begin{equation}
   \left( \bm{K}  f \right) (t) := \int_{0}^{T} K(t,s) f(s) \d s, \quad t \in [0,T], \quad f \in L^{2},
\end{equation}
and we denote by $\bm K^{\star}$ the dual operator  in $L^{2}$, that is the operator $\bm K^{\star}$ induced by the kernel $K^{\star}(t,s)=K(s,t)$, and by $\bm{id}$ the identity operator, i.e.~$(\bm{id} \,  u)(.) = u_{.}$.
The unique  solution $\hat{u} \in \mathcal L^2$ to the Fredholm equation \eqref{eq:Fredholm_equation} is then  explicitly given  by 
\begin{equation} \label{eq:explicit_optimal_control}
    \hat{u}_{t} = \left( \left( \bm{id} - \bm{B} \right)^{-1} a^{{\varphi},\alpha} \right)_{t}, \quad t \in [0,T],
\end{equation}
with
\begin{align} \label{eq:a_expression}
    a_{t}^{{\varphi},\alpha} & := \varphi_{t} + \alpha_{t} - \langle \mathbb{1}_{\{t<.\}} K(.,t), \bm{D}_{t}^{-1} \mathbb{1}_{\{t<.\}} \E_t \left[ \varphi_{.} + \alpha_{.} \right] \rangle_{L^2}, \\ \label{eq:B_expression}
    B(t,s) & := \mathbb{1}_{\{s<t\}} \left( \langle \mathbb{1}_{\{t<.\}} K(.,t), \bm{D}_{t}^{-1} \mathbb{1}_{\{t<.\}} K(.,s) \rangle_{L^2} - K(t,s) \right), \\ \label{eq:D_expression}
    \bm{D}_{t} & := \bm{id} + \bm{K}_{t} + (\bm{K}_{t})^{*},
\end{align}
 for $t,s \in [0,T]$. The explicit formula \eqref{eq:explicit_optimal_control} is derived in \cite{jaber2023equilibrium}[Proposition 5.1] in two steps:
 \begin{itemize}
     \item[(i)] for each fixed $t \in [0,T]$, the process $m_{t}(s) := \mathbb{1}_{\{ t \leq s\}} \E_{t} [u_{s}], \; s \in [0,T]$ is obtained explicitly in terms of the process $u$ as a solution to a Fredholm equation of the second kind\footnote{Recall a Fredholm equation of the second kind is of the form $u_{t} - \int_{0}^{t} L(t,s) u_{s} ds = f_{t}, \; t \in [0,T]$ and admits the unique solution $u_{t} = \left( \left( \bm{id} - \bm{L} \right)^{-1} f \right)_{t}$ under some assumptions on the source term $f$ and the kernel $L$.},
     
     \item[(ii)] then injecting such solution for $m_{t}(.)$ into the initial Fredholm equation \eqref{eq:Fredholm_equation} and re-arranging the terms yields another Fredholm equation of the second kind on $u$ only, which is solved explicitly by \eqref{eq:explicit_optimal_control}.
 \end{itemize}

We are now ready to state our main theorems on the existence and uniqueness of the optimal control to \eqref{eq:problem_formulation}, together with its explicit formula \eqref{eq:explicit_optimal_control} in terms of Lagrange multipliers as sufficient and necessary conditions of optimality.

\begin{theorem} \label{T:generalized_kkt_theorem_necessary}
Assume $\mathcal{U} \neq \emptyset$ and let $K \in \mathcal K$ be an admissible kernel. 
Then, there exists a unique admissible control $\hat u \in \mathcal U$ such that 
\begin{equation}\label{eq:optimalJu}
    {\mathcal{J}}(\hat{u}) = \inf_{u \in \mathcal{U}} {\mathcal{J}}(u).
\end{equation}
If in addition $int(\mathcal{U}) \neq \emptyset$, then there exist Lagrange multipliers $\left(\lambda^{i}\right)_{i \in \{1,2,3,4\}} \in \left( \mathcal{L}_{+}^{1} \right)^{4}$ satisfying the decomposition \eqref{eq:aggregated_lagrange_multiplier} with $\varphi \in \mathcal{L}^{2}$, as well as the complementary slackness equations \eqref{eq:KKT_slackness_conditions}, and such that the optimal control $\hat u$  
is given explicitly by \eqref{eq:explicit_optimal_control}.
\end{theorem}

\begin{proof}
    See Section \ref{ss:proof_theorem_necessary_sufficient_condition}. 
\end{proof}

A notable feature of Theorem~\ref{T:generalized_kkt_theorem_necessary} is that the non-necessarily Markovian optimization problem \eqref{eq:problem_formulation} with a general kernel $K$ admits an explicit form in terms of Lagrange multipliers in $\mathcal L^1$ via the formula \eqref{eq:explicit_optimal_control}. However, although the theorem guarantees the existence of such Lagrange multipliers as a necessary condition of optimality, the proof is nonconstructive for $\hat{u}$.\\

Therefore, we show in the following theorem that it is indeed sufficient to construct  Lagrange multipliers with no singular parts to obtain the optimal control $\hat{u}$. Then we will elaborate on the numerical construction of these Lagrange multipliers: we provide one example where they can be explicitly computed, and present later on in Section~\ref{S:numerics} a general algorithm to obtain them in a stochastic setting.

\begin{theorem} \label{T:generalized_kkt_theorem_sufficient}
Let $K \in \mathcal K$ be an admissible kernel. Assume there exist Lagrange multipliers $\lambda := \left(\lambda^{i}\right)_{i \in \{1,2,3,4\}} \in \left( \mathcal{L}_{+}^{1} \right)^{4}$ such that the process $\varphi$ defined by   \eqref{eq:aggregated_lagrange_multiplier} is in $\mathcal{L}^{2}$.
Let $u^{*} \in \mathcal{L}^{2}$ be the solution to the stochastic Fredholm equation 
\begin{equation} \label{eq:assume_fredholm}
    u_{t}^{*} + \int_0^{t} K(t,s) u_{s}^{*} \d s + \int_t^T K(s,t) \E_{t} u_{s}^{*} \d s = \varphi_{t} + \alpha_{t} \quad \left( \d t \otimes \mathbb{P} \right) - a.e.
\end{equation}
If $u^{*} \in \mathcal U$ and $(u^{*},\lambda)$ satisfies the complementary slackness conditions \eqref{eq:KKT_slackness_conditions}, then,  $u^{*} = \hat{u}$, where $\hat u$ is the unique optimal control of Theorem~\ref{T:generalized_kkt_theorem_necessary}.
\end{theorem}

\begin{proof}
    See Section \ref{ss:proof_sufficient_condition}.
\end{proof}

A key strength of Theorem \ref{T:generalized_kkt_theorem_sufficient} is that it is enough to construct four non-negative integrable stochastic Lagrange processes satisfying the slackness conditions \eqref{eq:KKT_slackness_conditions} and the decomposition \eqref{eq:aggregated_lagrange_multiplier} and such that the solution $u^*$ \eqref{eq:explicit_optimal_control} lies in $\mathcal{U}$, in order to obtain the unique optimal control, without the need for the Slater condition $int(\mathcal U)\neq \emptyset$. In one special case, we can actually derive explicitly the Lagrange multipliers as shown in the following example.

 \begin{example}[Explicit Lagrange multipliers] \label{ex:explicit_Lagrange_multipliers} 
Assume there is no transient impact ($K=0$) and only the trading rate $u$ is constrained with $a^{1}\leq u \leq a^{2}$. Indeed, taking $a^3$  small enough and $a^4$ large enough, then it is straightforward to check that the measures
    \begin{align}\label{eq:psiexplicit1}
        \psi^{1} (\d t,\d \omega) \;  = \;  \left( a_{t}^{1} - \alpha_{t} \right)^{+} \d t \P(\d \omega), \quad 
        \psi^{2} (\d t,\d \omega) \;  = \; \left( \alpha_{t} - a_{t}^{2} \right)^{+} \d t \P(\d \omega),\;  \psi^3=\psi^4\equiv 0,
    \end{align}
    and the control $u^*:=\alpha + \varphi $ with $\varphi:= \left( a^{1} - \alpha \right)^{+} - \left( \alpha-a^2 \right)^{+}  $ satisfy the conditions of Theorem~\ref{T:generalized_kkt_theorem_sufficient}, so that an application of the theorem yields that $\hat u=u^*$.  To be more precise,  constraining $u$ by functions $a^1$ and $a^2$ ensures that the inventory $X^u$  is  constrained  by $X_0 +\int_0^. a^1_s \d s \leq X^u_.\leq X_0 +\int_0^. a^2_s \d s $ so that  Theorem~\ref{T:generalized_kkt_theorem_sufficient} is applied for any $(a^3,a^4)$ such that $a^3_.< X_0 +\int_0^. a^1_s \d s$ and $a^4_. > X_0 +\int_0^. a^2_s \d s$.

\end{example}

In the general case however, we do not expect to provide such explicit formulas for the Lagrange multipliers (and hence for the unique optimal control given by equation \eqref{eq:explicit_optimal_control}), and we propose in Section~\ref{S:numerics} a stochastic version of Uzawa's algorithm for computing admissible sample trajectories of the optimal control by constructing step-by-step those of the Lagrange multipliers via a dual-primal scheme.

\begin{remark}[KKT optimality conditions and strong duality]
    \noindent Theorems~\ref{T:generalized_kkt_theorem_necessary} and \ref{T:generalized_kkt_theorem_sufficient} equivalently state that, under the Slater condition $int(\mathcal{U}) \neq \emptyset$, $\left( \hat{u}, \hat \lambda \right) \in \mathcal{L}^{\infty} \times (\mathcal L^1_+)^4$ satisfies the four KKT optimality conditions 
    \begin{itemize}
        \item[(i)] primal feasibility: $\hat{u} \in \mathcal{U}$,
        \item[(ii)] dual feasibility: $\hat \lambda^{i} \geq 0, \; i \in \{1,2,3,4\}$,
        \item[(iii)] slackness equations from \eqref{eq:KKT_slackness_conditions},
        \item[(iv)] and the Lagrangian stationary equation \eqref{eq:lagstation},
    \end{itemize}
    if and only if $\left( \hat{u}, \hat \lambda \right)$ is a saddle-point of the Lagrange functional $\mathcal{L}$ defined by \eqref{eq:def_lagrange_functional}, i.e.
    \begin{equation}
        \mathcal{L}\left( \hat{u}, \lambda \right) \leq \mathcal{L} \left(\hat{u}, \hat \lambda \right) \leq \mathcal{L}\left(v, \hat \lambda \right), \quad \lambda \in \left(\mathcal{L}_{+}^{1}\right)^{4}, \quad v \in \mathcal{U}.
    \end{equation}
\end{remark}

\section{ Numerical implementation} \label{S:numerics}

{In this section, we propose an algorithm to construct step-by-step admissible trajectories of the optimal control $\hat{u}$ and its associated optimal Lagrange multipliers from \eqref{eq:explicit_optimal_control} extending the standard Uzawa algorithm introduced in \cite{arrow1958studies} into the stochastic realm. Such algorithm is applied to solve various specifications of our stochastic control problem \eqref{eq:problem_formulation} in Section \ref{S:Applications}.

\subsection{Stochastic Uzawa algorithm}

We start by giving the general ideas behind the algorithm. Observe that the expression of the deterministic kernel $B$ in \eqref{eq:B_expression} depends only on the deterministic kernel $K$, while the stochastic process $a^{{\varphi},\alpha}$ in \eqref{eq:a_expression} depends on the kernel $K$, the stochastic signal $\alpha$ as well as the stochastic aggregated Lagrange multiplier $\varphi$ together with their conditional expectations. Recall that $\varphi$ from \eqref{eq:aggregated_lagrange_multiplier} is itself an explicit function of the Lagrange multipliers $\lambda := \left( \lambda^{i} \right)_{i \in \{ 1,2,3,4 \}}$ which account for the four inequality constraints in \eqref{eq:def_admissible_set_of_controls}.
\begin{itemize}
    \item[(i)] First, assume that $(\alpha, \varphi)$ together with their conditional expectations are known. Then, in the spirit of \cite[Section 5.1]{abi2022optimal}, a straightforward discretisation in time of \eqref{eq:explicit_optimal_control}, known as the Nyström method recalled in Section~\ref{s:solving_fredholm}, yields a discrete time approximation solution denoted by $\hat u^N \left( \lambda \right) =(\hat{u}_{t_{0}}, \hat{u}_{t_{1}}, \ldots, \hat{u}_{t_{N-1}})^\top$ on the uniformly spaced time grid
    \begin{equation} \label{eq:subdivision_def}
        \mathcal{T}_{N} := \left\{ 0 = t_{0} < t_{1} < \cdots < t_{N-1} < t_{N} = T \right\}, \quad N \in \mathbb{N}^{*},
    \end{equation}
    with step $\Delta := \frac{T}{N}$, such that 
    \begin{equation} \label{eq:explicit_discrete_time_control2}
        \hat{u}^{N}(\lambda)(\omega) \; = \; (I^{N}-B^{N})^{-1} a^{N, \alpha, \varphi}(\omega), \quad \omega \in \Omega,
    \end{equation}
    with $I^{N}$ the $N \times N$ identity matrix, $B^N$ a $N \times N$-explicit deterministic matrix and $a^{N, \alpha, \varphi}$ an $\mathbb R^{N}$-valued random variable given explicitly in terms of $(\alpha_{t_i}, \varphi_{t_i})_{0 \leq i \leq N-1}$ and $(\E_{t_i}[\alpha_{t_j}],\mathbb E_{t_i} [\varphi_{t_j}])_{0 \leq i < j \leq N}$, see Section \ref{s:solving_fredholm} for the precise expressions\footnote{Notice that when applying the discrete approximation to $\hat{u}$, we chose not to have a trading decision at maturity because it does not affect the final running inventory, nor the PnL functional \eqref{eq:stochasticConvexProgramFormulation}.}. Since the signal $\alpha$ is an input of the problem, it can be chosen in such a way that $(\E_{t_i}[\alpha_{t_j}])_{0\leq i < j\leq N}$ are explicitly known. However, the aggregated Lagrange multiplier $\varphi$ and its conditional expectations are in general not known and must also be approximated numerically on the grid $\mathcal{T}_{N}$. 

    \item[(ii)] Second, we construct step-by-step the unknown aggregated Lagrange multiplier $\varphi$ numerically using a stochastic adaptation of the Uzawa algorithm. Such an algorithm, described in detail in the next section, is simply an iterative, discrete time, stochastic gradient ascent applied on the dual formulation of the stochastic control problem:
    \begin{equation} \label{eq:dual_problem_formulation}
        \sup_{\lambda \in \mathcal (\mathcal L^1_+)^4} g(\lambda),
    \end{equation}
    where $g$ is the \textit{dual function} defined by
    \begin{equation} \label{eq:dual_function_def}
        g(\lambda) := \inf_{u \in \mathcal L^2} \mathcal{L} \left( u, \lambda \right), \quad \lambda \in \mathcal (\mathcal L^1)^4,
    \end{equation}
    where recall $\mathcal{L}$ is the Lagrangian defined by \eqref{eq:def_lagrange_functional}. We already know from Theorem \ref{T:generalized_kkt_theorem_necessary} that
    \begin{equation} \label{eq:minimizer_lagrangian_lbda_fixed}
        \hat{u}(\lambda) = \arg \inf_{u \in \mathcal L^2} \mathcal{L} \left( u, \lambda \right), \quad \lambda \in \mathcal (\mathcal L^1)^4,
    \end{equation}
    where $\hat{u}(\lambda)$ is explicitly given by \eqref{eq:explicit_optimal_control}. Consequently, we have
    \begin{equation} \label{eq:explicit_dual_function}
        g(\lambda) = \mathcal{L} \left( \hat{u}(\lambda), \lambda \right),
    \end{equation}
    and differentiating the Lagrangian functional \eqref{eq:def_lagrange_functional} with respect to its second variable yields the following projected ascent update when maximizing the dual function \eqref{eq:explicit_dual_function}
    \begin{equation} \label{eq:gradient_ascent_uzawa}
        \lambda \leftarrow \left( \lambda + \delta \mathcal{G}(\hat{u}(\lambda))\right)^{+},
    \end{equation}
    with $\delta>0$ some learning rate. Numerically, given some discrete time Lagrange multipliers $\lambda[n] = \lambda^{N}[n] := \left( \lambda_{t_{i}}^{1}[n], \lambda_{t_{i}}^{2}[n], \lambda_{t_{i+1}}^{3}[n], \lambda_{t_{i+1}}^{4}[n] \right)_{i \in \{ 0, \cdots, N-1 \}}$ at iteration $n \in \mathbb{N}^{*}$ of the gradient ascent\footnote{Notice we drop the dependence on the number of time steps $N$ to denote the discrete time Lagrange multipliers $\lambda[n]$ for the sake of clarity, and no ambiguity should arise given the context in what follows.}, the discrete time approximation $\varphi^{N}(\lambda[n])$ of the aggregated Lagrange multiplier $\varphi$ from \eqref{eq:aggregated_lagrange_multiplier} is given by
    \begin{equation} \label{eq:discrete_time_aggregated_lagrange_multiplier}
         \varphi_{t_{i}}^{N}(\lambda[n]) := \lambda_{t_{i}}^{1}[n] - \lambda_{t_{i}}^{2}[n] + \mathbb{E}_{t_{i}} \left[ \sum_{j=i}^{N-1} \left( \lambda_{t_{j+1}}^{3}[n] - \lambda_{t_{j+1}}^{4}[n] \right) \Delta \right], \quad i \in \{ 0, \cdots, N-1 \}.
    \end{equation}
    Notice that the time indices of the Lagrange multipliers are consistent with the facts that:
    \begin{enumerate}
        \item[a)] there is no trading decision to be constrained at maturity;
        \item[b)] the initial inventory satisfies the constraints at time $t=0$, recall \eqref{eq:def_constraining_processes}.
    \end{enumerate}
    Then, it follows from the approximation of $\hat{u}(\lambda)$ by \eqref{eq:explicit_discrete_time_control2} in (i) and the update equation \eqref{eq:gradient_ascent_uzawa} that stochastic Uzawa's discrete time formulation writes
    \begin{equation}
    \begin{cases} 
         \lambda[0] \; = \;  0, \\ 
         \hat{u}^{N}(\lambda[n]) \; = \; (I^{N}-B^{N})^{-1} a^{N}(\lambda[n]), \\
         \lambda[n+1] \; = \; \left( \lambda[n] + \delta_{n} \mathcal{G}^{N}(\hat{u}^{N}(\lambda[n])) \right)^{+},
     \end{cases}   
    \end{equation}
    where $\delta_n$ denotes the learning rate specified as\footnote{Such specification is inspired from the Robins-Monro convergence theorem initially tackled in \cite{robbins1951stochastic}.}
    \begin{equation} \label{eq:adaptive_learning_step}
        \delta_{n} = \frac{\delta}{n^{\beta}}, \quad n \in \mathbb{N}, \quad \delta, \beta > 0,
    \end{equation}
    and $\mathcal{G}^{N}(\hat{u}^{N}(\lambda[n]))$ is given by \begin{equation}
        \mathcal{G}^{N}(\hat{u}^{N}(\lambda[n])) := \begin{pmatrix}
          \left( a_{t_{i}}^{1} - \hat{u}_{t_{i}}^{N}(\lambda[n]) \right)_{i \in \{ 0, \cdots, N-1 \}} \\
         \left( \hat{u}_{t_{i}}^{N}(\lambda[n]) - a_{t_{i}}^{2} \right)_{i \in \{ 0, \cdots, N-1 \}} \\
         \left( a_{t_{i+1}}^{3} - X_{t_{i+1}}^{\hat{u}^{N}(\lambda[n])} \right)_{i \in \{ 0, \cdots, N-1 \}} \\
         \left( X_{t_{i+1}}^{\hat{u}^{N}(\lambda[n])} - a_{t_{i+1}}^{4} \right)_{i \in \{ 0, \cdots, N-1 \}}
         \end{pmatrix},
    \end{equation}
    with $\left(X_{t_{i}}^{\hat{u}^{N}(\lambda[n])}\right)_{0 \leq i \leq N}$ denoting the inventory associated to the discrete time control $\hat{u}^{N}(\lambda[n])$ such that
    \begin{equation} \label{eq:discrete_time_inventory}
        X_{t_{i}}^{\hat{u}^{N}(\lambda[n])} := X_{0} + \Delta \sum_{l=0}^{i-1} \hat{u}_{t_{l}}^{N}(\lambda[n]), \quad 0 \leq i \leq N.
    \end{equation}
    Notice that $\hat{u}^{N}(\lambda[n])$ is simply the analogue of \eqref{eq:explicit_discrete_time_control2} but with the input $\varphi^{N}(\lambda[n])$ which is now known at iteration $n$ thanks to \eqref{eq:discrete_time_aggregated_lagrange_multiplier}. Indeed, $a^{N}(\lambda[n]))$, given in \eqref{eq:source_term_discrete_time}, depends explicitly on $(\varphi_{t_i}^{N}(\lambda[n]))_{0 \leq i \leq N-1}$ and its conditional expectations $\left( \mathbb E_{t_{i}} \left[\varphi_{t_{j}}^{N}(\lambda[n]) \right] \right)_{0\leq i < j\leq N}$ which, in our applications, are approximated using least-squares Monte Carlo regressions in the same spirit as \citet{longstaff2001valuing}.
\end{itemize}
The pseudocode is detailed in Algorithm \ref{algo:stochatic_uzawa}.

\begin{algorithm}
\caption{Stochastic Uzawa}
\label{algo:stochatic_uzawa}
\small{
\textbf{Objective.} Construct admissible sample trajectories of the optimal control solution to \eqref{eq:problem_formulation} by solving step-by-step the discrete time approximation of the dual problem formulation \eqref{eq:dual_problem_formulation}.

\textbf{Inputs.} 
\begin{itemize}
    \item Fix $N \in \mathbb{N}^{*}$ the number of time steps in the grid $\mathcal{T}_{N}$ given by \eqref{eq:subdivision_def}, and $M \in \mathbb{N}^{*}$ the number of trajectories.
    \item An admissible kernel $K \in \mathcal{G}$.
    \item A function that can generate sample paths of the price signal $(\alpha_{t_{i}})_{0 \leq i \leq N}$, its conditional expectations $(\E_{t_{i}}[\alpha_{t_{j}}])_{0 \leq i < j \leq N}$ and sample paths of the constraining processes $((a_{t_{j}}^{i})_{0 \leq j \leq N})_{i \in \{1, 2, 3, 4 \}} \in \left(\mathcal{L}^{\infty}\right)^{4}$ satisfying \eqref{eq:def_constraining_processes} at the dates of $\mathcal{T}_{N}$.
    \item Fix a maximum number of iterations $D \in \mathbb{N}^{*}$, a tolerance threshold $\bar \epsilon > 0$ and some parameters $\delta, \beta > 0$ for the adaptive learning step given by \eqref{eq:adaptive_learning_step}.
\end{itemize}

\vspace{-0.1cm}
\textbf{Outputs.}
\begin{itemize}
    \item $M$ sample trajectories on $\mathcal{T}_{N}$ of the approximated optimal control $\hat{u}^{N}$ of \eqref{eq:explicit_optimal_control}, and those of the associated Lagrange multipliers $\lambda := \left(\lambda^{i}\right)_{i \in \{1,2,3,4\}}$.
\end{itemize}

\vspace{-0.1cm}
\textbf{Variable Declarations.}
\begin{itemize}
        \item Initialize an integer $n$ used as an iteration rank tracker.
        \item Initialize empirical slackness variables $\bar{S} := (\bar{S}^{i})_{i \in \{1, 2, 3, 4 \}} \in \mathbb{R}^{4}$ to quantify how well the slackness conditions \eqref{eq:KKT_slackness_conditions} are satisfied.
        \item Initialize four times $M$ trajectories $\lambda[n] := \left( \lambda^{1}, \lambda^{2}, \lambda^{3}, \lambda^{4} \right) \in \mathbb{R}^{4 \times M \times N}$ to store the sample trajectories of the four Lagrange multipliers at the dates of $t_{0}, \cdots, t_{N-1}$ for $\lambda^{1}, \lambda^{2}$ and at the dates $t_{1}, \cdots, t_{N}$ for $\lambda^{3}, \lambda^{4}$.
\end{itemize}

\vspace{-0.1cm}
\textbf{Algorithm steps.}
\begin{enumerate}
    \item Generate $M$ trajectories of the price signal $(\alpha_{t_{i}})_{0 \leq i \leq N}$, its conditional expectations $(\E_{t_{i}}[\alpha_{t_{j}}])_{0 \leq i < j \leq N}$, and $M$ trajectories of the constraining processes at the dates of $\mathcal{T}_{N}$.
    
    \item Initialize $n = 0$, $\bar{S}[0] = \left( \infty, \infty, \infty, \infty \right)$ and the Lagrange multipliers by $\lambda[0] = \left( 0, 0, 0, 0 \right)$.

    \item Approximate the minimizer of \eqref{eq:minimizer_lagrangian_lbda_fixed} by $\hat{u}^{N}(\lambda[0])$ using the Nyström scheme via \eqref{eq:explicit_discrete_time_control3}, and compute the associated inventory $\left(X_{t_{i+1}}^{\hat{u}^{N}(\lambda[0])}\right)_{0 \leq i \leq N-1}$ using \eqref{eq:discrete_time_inventory}.
    
    \item While $n < D$ and $\max \left( (\bar{S}^{i}[n])_{i \in \{1, 2, 3, 4 \}} \right) > \bar \epsilon$, do:
    \begin{enumerate}

        \item Update the Lagrange multipliers' sample paths in the direction of increasing gradient of the dual function \eqref{eq:dual_function_def} projected onto the positive numbers:
        \vspace{-0.15cm}
        \begin{equation} \label{eq:explicit_discrete_time_control_uzawa}
            \begin{pmatrix} 
         \lambda_{t_{i}}^{1}[n+1] \\
         \lambda_{t_{i}}^{2}[n+1] \\
         \lambda_{t_{i+1}}^{3}[n+1] \\
         \lambda_{t_{i+1}}^{4}[n+1]
         \end{pmatrix} \; = \;  \begin{pmatrix}
         \left( \lambda_{t_{i}}^{1}[n] + \delta_{n} \left( a_{t_{i}}^{1} - \hat{u}^{N}_{t_{i}}(\lambda[n]) \right) \right)^{+} \\
         \left( \lambda_{t_{i}}^{2}[n] + \delta_{n} \left( \hat{u}^{N}_{t_{i}}(\lambda[n]) - a_{t_{i}}^{2} \right) \right)^{+} \\
         \left( \lambda_{t_{i+1}}^{3}[n] + \delta_{n} \left( a_{t_{i+1}}^{3} - X_{t_{i+1}}^{\hat{u}^{N}(\lambda[n])} \right) \right)^{+}\\
         \left( \lambda_{t_{i+1}}^{4}[n] + \delta_{n} \left( X_{t_{i+1}}^{\hat{u}^{N}(\lambda[n])} - a_{t_{i+1}}^{4} \right) \right)^{+}
         \end{pmatrix}, \quad i \in \{0, \cdots, N-1\}.
        \end{equation}
        
        \item Estimate $\left( \varphi_{t_{i}}^{N}(\lambda[n+1]) \right)_{0 \leq i \leq N-1}$ and its conditional expectations $(\E_{t_{i}}[\varphi_{t_{j}}^{N}(\lambda[n+1])])_{0 \leq i < j \leq N}$ via least-squares Monte Carlo regressions using respectively \eqref{eq:E_j_linear_approx}--\eqref{eq:E_j_k_linear_approx}, plugging $\lambda[n+1]$ in the expression of $\varphi^{N}(\lambda[n+1])$ given by \eqref{eq:discrete_time_aggregated_lagrange_multiplier}.
        
       \item Approximate the minimizer of \eqref{eq:minimizer_lagrangian_lbda_fixed} by $\hat{u}^{N}(\lambda[n+1])$ using the Nyström scheme via \eqref{eq:explicit_discrete_time_control3}, and compute the associated inventory $\left(X_{t_{i+1}}^{\hat{u}^{N}(\lambda[n+1])}\right)_{0 \leq i \leq N-1}$ using \eqref{eq:discrete_time_inventory}.
    
        \item Compute the discrete time empirical slackness variables
        \vspace{-0.15cm}
        \begin{equation} \label{eq:discrete_time_complementary_slackness}
            \begin{cases}
                \bar{S}^{1}[n+1] & = \frac{1}{M} \sum_{m = 1}^{M} \left[ \left( a_{t_{i}}^{1} \left( \omega_{m} \right) - \hat{u}^{N}_{t_{i}}(\lambda[n+1]) \left( \omega_{m} \right) \right) \lambda_{t_{i}}^{1} \left( \omega_{m} \right) \right] \\
                \bar{S}^{2}[n+1] & = \frac{1}{M} \sum_{m = 1}^{M} \left[ \left( \hat{u}^{N}_{t_{i}}(\lambda[n+1]) \left( \omega_{m} \right) - a_{t_{i}}^{2} \right) \lambda_{t_{i}}^{2} \left( \omega_{m} \right) \right] \\
                \bar{S}^{3}[n+1] & = \frac{1}{M} \sum_{m = 1}^{M} \left[ \left( a_{t_{i+1}}^{3} \left( \omega_{m} \right) - X_{t_{i+1}}^{\hat{u}^{N}(\lambda[n+1])} \left( \omega_{m} \right) \right) \lambda_{t_{i+1}}^{3} \left( \omega_{m} \right) \right] \\
                \bar{S}^{4}[n+1] & = \frac{1}{M} \sum_{m = 1}^{M} \left[ \left( X_{t_{i+1}}^{\hat{u}^{N}(\lambda[n+1])}\left( \omega_{m} \right) - a_{t_{i+1}}^{4}\left( \omega_{m} \right) \right) \lambda_{t_{i+1}}^{4}\left( \omega_{m} \right) \right]
            \end{cases}, \quad i \in \{ 0, \cdots, N-1 \},
        \end{equation}
        where $\omega_{m} \in \Omega, \; m \in \{1, \cdots, M\}$ is the $m^{th}$ random element associated to the $m^{th}$ sample trajectory of $x$.

        \item Increment the rank tracker $n$ of one unit.
    \end{enumerate}
    
    \item Return the $M$ sample paths of $\left( \hat{u}^{N}(\lambda[n]), \lambda[n] \right)$.
\end{enumerate}
}
\end{algorithm}

\subsection{Solving the Fredholm equation via the Nyström scheme: steps 3. and 4.c) of Algorithm \ref{algo:stochatic_uzawa}} \label{s:solving_fredholm}

The Nyström approximation scheme consists in approximating the action of linear operators by matrix operations.  The key underlying approximation step relies on a left-rectangle integral approximation such that the solution $u \in \mathcal{L}^{2}$ to the following Fredholm equation
\begin{equation}
    u_{t} + \int_{0}^{t} B(t,s) u_{s} \d s = f_{t}, \quad t \in [0,T],
\end{equation}
for some admissible kernel $B$ (e.g.~the one given by \eqref{eq:B_expression}), and some square-integrable process $f$ (e.g.~$a^{\alpha, \varphi}$ from \eqref{eq:a_expression}), is approximated on the time grid $\mathcal{T}_{N}$ \eqref{eq:subdivision_def} by
\begin{equation}
    u_{t_{i}} + \sum_{j=0}^{i-1} \left( \int_{t_{j}}^{t_{j+1}} B(t_{i},s) \d s \right) u_{t_{j}} = f_{t_{i}}, \quad i \in \{0, \cdots, N-1\}.
\end{equation}
Consequently, we can write
\begin{equation}
    u^{N} = \left( I^{N} - B^{N} \right)^{-1} f^{N},
\end{equation}
where $u^{N} := (u_{t_{0}}, u_{t_{1}}, \ldots, u_{t_{N-1}})^\top$, $f^{N} := (f_{t_{0}}, f_{t_{1}}, \ldots, f_{t_{N-1}})^\top$, $I^{N}$ is the $N \times N$ identity matrix and $B^{N}$ is the lower-triangular matrix\footnote{Notice the matrix $B^{N}$ is formed by the respective integrals of the kernel $B$ on the different intervals of $\mathcal{T}_{N}$, which makes the scheme numerically stable even for integrable kernels that blow up at the origin, such as the fractional kernel \eqref{eq:def_fractional_kernel}.} given by \begin{equation}
    B_{i,j}^{N} := \mathbb{1}_{0 \leq j < i \leq N-1} \int_{t_{j}}^{t_{j+1}} B(t_{i},s) \d s, \quad i, j = 0, \ldots, N-1.
\end{equation}
The precise expressions of $B$ and $f$ in our context are given in what follows. For this, we introduce matrices related to the transient impact kernel $K$:
\begin{align}
    L_{i,j}^{N} & := \mathbb{1}_{0 \leq j < i \leq N-1} \int_{t_j}^{t_{j+1}} K(t_i,s) \d s, \quad i, j = 0, \ldots, N-1, \\
    U_{i,j}^{N} & := \mathbb{1}_{0 \leq i \leq j \leq N-1} \int_{t_j}^{t_{j+1}} K(s,t_i) \d s, \quad i, j = 0, \ldots, N-1, \\
    L_{k,j}^{N,i} & := L_{k,j}^{N} 1_{k \leq i \leq (N-1)} 1_{j \leq i \leq (N-1)}, \quad k,j = 0,\ldots, N-1, \\
    U_{k,j}^{N,i} & := U_{k,j}^{N} 1_{k \leq i \leq (N-1)} 1_{j \leq i \leq (N-1)}, \quad k,j = 0,\ldots, N-1.
\end{align}
Then, in steps 3. and 4.c) of Algorithm \ref{algo:stochatic_uzawa}, we approximate the optimal solution \eqref{eq:explicit_optimal_control} by the expression \eqref{eq:explicit_discrete_time_control2} which, in the context of the algorithm at step $n \in \mathbb{N}$, writes as
\begin{equation} \label{eq:explicit_discrete_time_control3}
    \hat{u}_{t_{i}}^{N}(\lambda[n+1]) \; = \; \left( (I^{N}-B^{N})^{-1} a^{N, \alpha, \varphi^{N}(\lambda[n+1])} \right)_{i}, \quad i \in \{0, \cdots, N-1\}, 
\end{equation}
where we defined
\begin{align} \label{eq:source_term_discrete_time}
\hspace{-0.7cm}
    a^{N, \alpha, \varphi^{N}(\lambda[n+1])} & := \left( R_{i,i}^{N,\alpha} + \Lambda_{i,i}^{N}(\lambda[n+1]) - \left(U^{N}\right)^{i.} (D^{N,i})^{-1} \left( R^{N,\alpha} + \Lambda^{N}(\lambda[n+1]) \right)^{.i} \right)_{i \in \{0, \ldots, N-1\}}, \\
    \left(B^{N}\right)_{i,j} & := \left(U^{N}\right)^{i.} (D^{N,i})^{-1} \left(L^{N}\right)^{.j} - L_{ij}^{N}, \quad i,j=0,\ldots, N-1, \\
    D^{N,i} & := I^{N} + L^{N,i} + U^{N,i}, \quad i=0,\ldots, N-1,
\end{align}
where $X^{i.}$ denotes the $i^{th}$ vector-row of $X$ and $X^{.j}$ its $j^{th}$ vector-column, and where $R^{N,\alpha}$ is a triangular matrix containing the alpha process $\alpha$ and its conditional expectations $(\E_{t_{i}}[\alpha_{t_{j}}])_{0 \leq i < j \leq N}$ such that
\begin{equation} \label{eq:cond_exp_signal}
    R_{i,j}^{N,\alpha} := \mathbb{1}_{0 \leq i \leq j \leq (N-1)} \mathbb{E}_{t_{i}} \left[ \alpha_{t_{j}} \right], \quad i,j=0,\ldots, N-1,
\end{equation}
while $\Lambda^{N}(\lambda[n+1])$ is similarly given by
\begin{align}
    \Lambda_{i,j}^{N}(\lambda[n+1]) & := \mathbb{1}_{0 \leq i \leq j \leq (N-1)} \mathbb{E}_{t_{i}} \left[ \varphi_{t_{j}}^{N}(\lambda[n+1]) \right] \\ \label{eq:cond_multipliers_discrete}
    & = 1_{i \leq j \leq (N-1)} \mathbb{E}_{t_{i}} \left[ \lambda_{t_{j}}^{1}[n+1] - \lambda_{t_{j}}^{2}[n+1] + \sum_{l=j+1}^{N}  \left( \lambda_{t_{l}}^{3}[n+1] - \lambda_{t_{l}}^{4}[n+1] \right) \right], \quad i,j = 0,\ldots, N-1,
\end{align}
with the second equality stemming from \eqref{eq:discrete_time_aggregated_lagrange_multiplier}.

\subsection{Estimation of the conditional expectations in $\Lambda^{N}(\lambda[n+1])$: step 4.b) of Algorithm \ref{algo:stochatic_uzawa}} \label{ss:estimation_conditional_expectations}

Notice that in order to compute the optimal control $\hat{u}^{N}(\lambda[n+1])$ given by \eqref{eq:explicit_discrete_time_control3} in step 4.b) of Algorithm \ref{algo:stochatic_uzawa}, we first have to estimate $\left( \varphi_{t_{i}}^{N}(\lambda[n+1]) \right)_{0 \leq i \leq N-1}$ from \eqref{eq:discrete_time_aggregated_lagrange_multiplier} and its conditional expectations $(\E_{t_{i}}[\varphi_{t_{j}}^{N}(\lambda[n+1])])_{0 \leq i < j \leq N}$. Noting that only expectations are to be computed for $i=0$ in \eqref{eq:cond_multipliers_discrete} and that $\lambda_{t_{i}}^{1}, \lambda_{t_{i}}^{2}$ are $\mathcal{F}_{t_{i}}$-measurable for any $0 < i \leq (N-1)$, then we only have to estimate the respective conditional expectations
\begin{align}
    \label{eq:E_j}
    E_{i}[n+1] & := \mathbb{E}_{t_{i}} \left[ \sum_{l=i+1}^{N} \left( \lambda_{l}^{3}[n+1] - \lambda_{l}^{4}[n+1] \right) \right], \quad 0 < i \leq (N-1),\\ \label{eq:E_j_k}
    \hspace{-4cm}
    E_{i,j}[n+1] & := \mathbb{E}_{t_{i}} \left[ \lambda_{j}^{1}[n+1] - \lambda_{j}^{2}[n+1] + \sum_{l=j+1}^{N}  \left( \lambda_{l}^{3}[n+1] - \lambda_{l}^{4}[n+1] \right) \right], \quad 0 < i < j \leq (N-1),
\end{align}
to obtain respectively the diagonal and off-diagonal terms of $\Lambda^{N}(\lambda[n+1])$. Least-square Monte Carlo regressions can be used once the Markovian variables of the problem have been identified. For instance, if we assume the constraining functions to be deterministic, the signal to be Markovian and the transient impact kernel $K$ to be exponential as in \eqref{eq:def_exponential_kernel}, then the underlying randomness comes only from the stochastic signal, and only one time-scale has to be taken into account for the impact decay. On the other hand, if the constraints turn out to be stochastic, and if a power law or multiple time scales are considered in the impact decay, then the regression variables should be adapted accordingly\footnote{For example, the power law impact decay can be approximated by multiple time scales, see \cite{abi2019multifactor} for example.}. In this paper, based on the Lagrange multipliers' update in \eqref{eq:explicit_discrete_time_control_uzawa}, we chose\footnote{In order to approximate the filtration $\mathcal{F}_{t_{i}}$ in $\mathbb{E}_{t_{i}}$, we tested various combinations among the dependent variables $\left( \alpha_{t_{i}}, Z_{t_{i}}^{\hat{u}^{N}(\lambda[n])}, X_{t_{i}}^{\hat{u}^{N}(\lambda[n])}, \hat{u}_{t_{i}}^{N} \right)$ and found out the most stable numerical results were observed for the choice made in \eqref{eq:E_j_approx}--\eqref{eq:E_j_k_approx}.} to approximate
\begin{align} \label{eq:E_j_approx}
    E_{i}[n+1] & \approx \chi_{i}^{n+1}(\alpha_{t_{i}}, Z_{t_{i}}^{\hat{u}^{N}(\lambda[n])}, X_{t_{i}}^{\hat{u}^{N}(\lambda[n])}), \quad 0 < i \leq (N-1), \\
    \label{eq:E_j_k_approx}
    E_{i,j}[n+1] & \approx \phi_{i,j}^{n+1}(\alpha_{t_{i}}, Z_{t_{i}}^{\hat{u}^{N}(\lambda[n])}, X_{t_{i}}^{\hat{u}^{N}(\lambda[n])}), \quad 0 < i < j \leq (N-1),
\end{align}
where $\left( \phi_{i,j}^{n+1} \right)_{0 < i < j \leq (N-1)}$ and $\left( \chi_{i}^{n+1} \right)_{0 < i \leq (N-1)}$ are deterministic measurable functions of the signal $\left( \alpha_{t_{i}} \right)_{0 < i \leq (N-1)}$, the Markovian state variable 
\begin{equation}
    Z_{t_{i}}^{\hat{u}^{N}(\lambda[n])} := \langle L^{i.}, \hat{u}^{N}(\lambda[n]) \rangle, \quad 0 < i \leq (N-1),
\end{equation}
which is the left-rectangle approximation of the running transient impact $Z^{u}$ from \eqref{eq:transient_market_impact_state_variable}, and the running inventory $\left( X_{t_{i}}^{\hat{u}^{N}(\lambda[n])} \right)_{0 < i \leq (N-1)}$, in the spirit of Doob's Theorem \cite[35, Chapter 1, p. 18]{dellacherie1975probabilit}. Notice that the last two variables are built using $\hat{u}^{N}(\lambda[n])$ obtained from the previous step of the algorithm. More precisely, we specify respectively the functions $\left( \phi_{i,j}^{n+1} \right)_{0 < i < j \leq (N-1)}$ and $\left( \chi_{i}^{n+1} \right)_{0 < i \leq (N-1)}$ in \eqref{eq:E_j_k_approx}--\eqref{eq:E_j_approx} as polynomial functions: for $d \in \mathbb{N}$, we fit the respective linear functions $\langle l_{i,j}^{\phi, d}[n+1], . \rangle$ and $\langle l_{i}^{\chi, d}[n+1], . \rangle$ such that real coefficients $\left( l_{i,j}^{\phi, d}[n+1] \right)$ and $\left( l_{i}^{\chi, d}[n+1] \right)$ are respectively obtained by minimizing the least-squares of the respective dependent variables
\begin{align}
    Y_{i}^{\chi}[n+1] & := \sum_{l=i+1}^{N} \left( \lambda_{t_{l}}^{3}[n+1] - \lambda_{t_{l}}^{4}[n+1] \right), \quad 0 < i \leq (N-1),\\
    Y_{j}^{\phi}[n+1] & := \lambda_{t_{j}}^{1}[n+1] - \lambda_{t_{j}}^{2}[n+1] + \sum_{l=j+1}^{N}  \left( \lambda_{t_{l}}^{3}[n+1] - \lambda_{t_{l}}^{4}[n+1] \right), \quad 0 < j \leq (N-1),
\end{align}
against a Laguerre polynomial basis expansion\footnote{As in \cite{longstaff2001valuing}, note that other basis expansions may well have been considered, for example Hermite, Legendre, Chebyshev, Gegenbauer, and Jacobi polynomials.} of the state variables - signal, running transient impact and running inventory - appearing in \eqref{eq:E_j_approx}--\eqref{eq:E_j_k_approx} with maximum degree $d$, defined for any $0 < i \leq (N-1)$ as the family
\begin{equation} \label{eq:basis_expansion}
    \mathcal{B}^d(\alpha_{t_{i}}, Z_{t_{i}}^{\hat{u}^{N}(\lambda[n])}, X_{t_{i}}^{\hat{u}^{N}(\lambda[n])}) := \left( L_{p} \left( \alpha_{t_{i}} \right) L_{q} \left( Z_{t_{i}}^{\hat{u}^{N}(\lambda[n])} \right) L_{r} \left( X_{t_{i}}^{\hat{u}^{N}(\lambda[n])} \right) \right)_{\underset{p, q, r \geq 0}{p+q+r \leq d}},
\end{equation} 
where $L_{p}$ denotes the standard Laguerre polynomial of degree $p \in \mathbb N$. Finally, we estimate the conditional expectations \eqref{eq:E_j} and \eqref{eq:E_j_k} by 
\begin{align} \label{eq:E_j_linear_approx}
    E_{i}[n+1] \approx \chi_{i}^{n+1}(\alpha_{t_{i}}, Z_{t_{i}}^{\hat{u}^{N}(\lambda[n])}, X_{t_{i}}^{\hat{u}^{N}(\lambda[n])}) & \approx \langle l_{i}^{\chi, d}[n+1], \mathcal{B}^d(\alpha_{t_{i}}, Z_{t_{i}}^{\hat{u}^{N}(\lambda[n])}, X_{t_{i}}^{\hat{u}^{N}(\lambda[n])}) \rangle, \quad 0 < i \leq (N-1), \\
    \label{eq:E_j_k_linear_approx}
    E_{i,j}[n+1] \approx \phi_{i,j}^{n+1}(\alpha_{t_{i}}, Z_{t_{i}}^{\hat{u}^{N}(\lambda[n])}, X_{t_{i}}^{\hat{u}^{N}(\lambda[n])}) & \approx \langle l_{i,j}^{\phi, d}[n+1], \mathcal{B}^d(\alpha_{t_{i}}, Z_{t_{i}}^{\hat{u}^{N}(\lambda[n])}, X_{t_{i}}^{\hat{u}^{N}(\lambda[n])}) \rangle, \quad 0 < i < j \leq (N-1),
\end{align}
using the $M \in \mathbb{N}^{*}$ sample trajectories of the discrete time signal process $\alpha$, as well as those of $\hat{u}^{N}(\lambda[n])$, computed at the previous step of the algorithm. Note that the dimension of the optimal regression coefficients $\left( l_{i,j}^{\phi, d}[n+1] \right)_{\{ 0 < i < j \leq (N-1) \}}$ and $\left( l_{i}^{\chi, d}[n+1] \right)_{\{ 0 < i \leq (N-1) \}}$ is given by the number of terms in the basis expansion family \eqref{eq:basis_expansion}.
}

\section{Application to optimal execution and optimal energy storage} \label{S:Applications}
In this section, we apply the stochastic Uzawa Algorithm \ref{algo:stochatic_uzawa} to the problems of optimal execution and energy storage. For this, we consider two independent Brownian motions $B$ and $W$, and we specify the exogenous price process $S$ from \eqref{eq:exogenous_price_process} as
\begin{equation} \label{eq:drift_like_signal_spec}
    S_{t} = S_{0} + \int_{0}^{t} I_{s} ds + \sigma B_{t},
\end{equation}
where $\sigma > 0$, and $I$ is a level-seasonal Ornstein-Uhlenbeck process with dynamics
\begin{equation} \label{eq:drift_signal_specification}
\d I_t = (A(t) - \kappa I_t) \d t + \xi \d W_t,
\end{equation}
with
\begin{equation} \label{eq:seasonal_level_specification}
A(t) := \theta \sin(wt+\phi), \quad \theta, w \geq 0, \quad \phi \in [0, 2 \pi).
\end{equation}
The signal-related conditional expectations matrix from \eqref{eq:cond_exp_signal} is explicitly given by 
\begin{align}
    R_{i,j}^{N, \alpha} \; = \; & 1_{0 \leq i \leq j \leq (N-1)} \left( \left( I_{t_{i}} - \frac{\theta}{\kappa^{2}+w^{2}} \left[ \kappa \sin(w t_i + \phi) - w \cos(w t_i + \phi) \right] \right) \frac{e^{- \kappa (t_{j}-t_{i})} - e^{- \kappa (T-t_{i})}}{\kappa} \right. \\
    & \left. - \frac{\theta}{\kappa^{2}+w^{2}} \left( \frac{\kappa}{w} \left(\cos(wT + \phi) - \cos(w t_{j} + \phi) \right) + \left( \sin(w T + \varphi) - \sin(wt_{j} + \phi) \right) \right) \right).
\end{align}
We also consider both the exponential and fractional Volterra kernels $K$ from \eqref{eq:def_exponential_kernel}--\eqref{eq:def_fractional_kernel} for the transient impact component. Finally, the time grid parameters are set as $T=1$ and $N=100$ in \eqref{eq:subdivision_def}, and the maximum polynomial degree used in for the Laguerre basis expansion \eqref{eq:basis_expansion} for the approximation of conditional expectations from \eqref{eq:E_j_k_linear_approx}--\eqref{eq:E_j_linear_approx} is set to $d=2$ for all the following numerical applications.

\subsection{Optimal execution under (stochastic) constraints} \label{ss:optimal_execution_under_constraints}

We first consider an optimal execution problem with transient impact, stochastic signal, with a full terminal liquidation constraint and under possible additional constraints on both the trading rate and the inventory. The problem takes the form:
\begin{align} \label{eq:optimal_liquidation}
\sup_{u \in \mathcal{U}} & \; \mathbb{E} \left[ \int_{0}^{T} \left( \alpha_{t} - \left( \frac{1}{2} u_{t} + \int_{0}^{t} K(t,s) u_{s} \d s \right) \right) u_{t} \d t \right], \\ \label{eq:optimal_liquidation_admissible_set}
& \mathcal{U} := \left\{ u \in \mathcal{L}^{\infty}, \; \; a^{1} \leq u \leq a^{2} \; \text{and} \; a^{3} \leq X^{u} \leq a^{4} \; \left( \d t \otimes \mathbb{P} \right) -a.e. \right\},
\end{align}
where $a^{1}, a^{2}, a^{3}, a^{4}$ are specified below depending on the four following configurations that we consider:
\begin{itemize}
    \item[(i)] a numerical sanity check for our numerical algorithm with no-additional constraints except for the terminal liquidation\footnote{Equality constraints are indeed permissible when building Lagrange multipliers, see Remark \ref{R:Slater_and_numerics}.};
    \item[(ii)] a `no-buying' constraint with terminal liquidation in the presence of a `buy' signal; 
    \item[(iii)] a `no-shorting' constraint with terminal liquidation in the presence of a `sell' signal;
    \item[(iv)] a `stop-trading' stochastic constraint whenever the exogenous price drops below a specified reference level.
\end{itemize}

\paragraph{Mimic the absence of constraints.} In what follows, for some constraints' specifications, we set large enough constants $\bar{M}, \bar{M}' >0$ such that the sample trajectories of the associated Lagrange multipliers are observed to be zeros (i.e.non binding constraints), except at maturity for the Lagrange multipliers $\lambda^{3}$ and $\lambda^{4}$ to enforce full liquidation. For all our numerical illustrations, a value around $1e16$ safely achieves this. Thus, the role of these large constants is to mimic numerically in our algorithm the case where the trading rate $u$ on $[0,T]$ and/or the inventory $X^u$ on $[0,T)$ are upper and/or lower unconstrained.


\begin{enumerate}
    \item[(i)] \textbf{Sanity check of the stochastic Uzawa Algorithm \ref{algo:stochatic_uzawa}.} First, we illustrate the convergence of our stochastic Uzawa algorithm when the agent only takes into account the full liquidation constraint by setting in \eqref{eq:optimal_liquidation_admissible_set}
    \begin{equation} \label{eq:assanity}
         -a_{t}^{1} = a_{t}^{2} = \bar M,  \quad t \in [0,T], \quad  -a_{t}^{3} =  a_{t}^{4} = \bar M, \quad t \in [0,T),\quad 
         a_{T}^{3} = a_{T}^{4} = 0.
    \end{equation}
    where recall $\bar M$ is a sufficiently large to mimic the absence of constraints except at the terminal date. Inspired by the problem formulation \cite[Equation (2.6)]{abi2022optimal}, we introduce as a benchmark the liquidation problem with soft quadratic terminal constraint in the form
    \begin{equation} \label{eq:benchmark_soft_constraint}
        \sup_{u \in \mathcal{L}^{2}} \mathbb{E} \left[ \int_{0}^{T} \left( \alpha_{t} - \left( \frac{1}{2} u_{t} + \int_{0}^{t} K(t,s) u_{s} \d s \right) \right) u_{t} \d t - \nu \left( X_{T}^{u} \right)^{2} \right],
    \end{equation}
    for some $\nu \geq 0$. The optimum of \eqref{eq:benchmark_soft_constraint} is known in explicit form, see \cite[Proposition 4.5]{abi2022optimal}. In Figure \ref{F:convergence_plot}, we compare the resulting control of \eqref{eq:benchmark_soft_constraint} when $\nu >> 1$, to enforce an approximate terminal liquidation constraint, to the one computed by our stochastic Uzawa algorithm for the problem \eqref{eq:optimal_liquidation} with the constraints \eqref{eq:assanity} for one sample trajectory. Notice that the empirical slackness variables $(\bar{S}^{i})_{i \in \{3, 4 \}}$ defined in \eqref{eq:discrete_time_complementary_slackness} converge fast to zero in the bottom-left plot ($(\bar{S}^{i})_{i \in \{1, 2\}}$ are zeros as the rate constraints never get activated, and are not plotted), while all the resulting sample controls satisfy the full liquidation constraint on all trajectories up to a precision of less 1e-6 after only $100$ iterations of the algorithm as illustrated by the bottom-right plot. Also, the numerical convergence in terms of violation gap is faster for the exponential impact decay than for the power law kernel, which could be explained by the choice of the regression variables used to estimate the conditional expectations from \eqref{eq:E_j_k_linear_approx}--\eqref{eq:E_j_linear_approx}, see Section \ref{ss:estimation_conditional_expectations}.
    
    \begin{figure}[H]
    \begin{center}
    \vspace*{-0.15in}
    \hspace*{-0.3in}
    \includegraphics[width=7 in,angle=0]{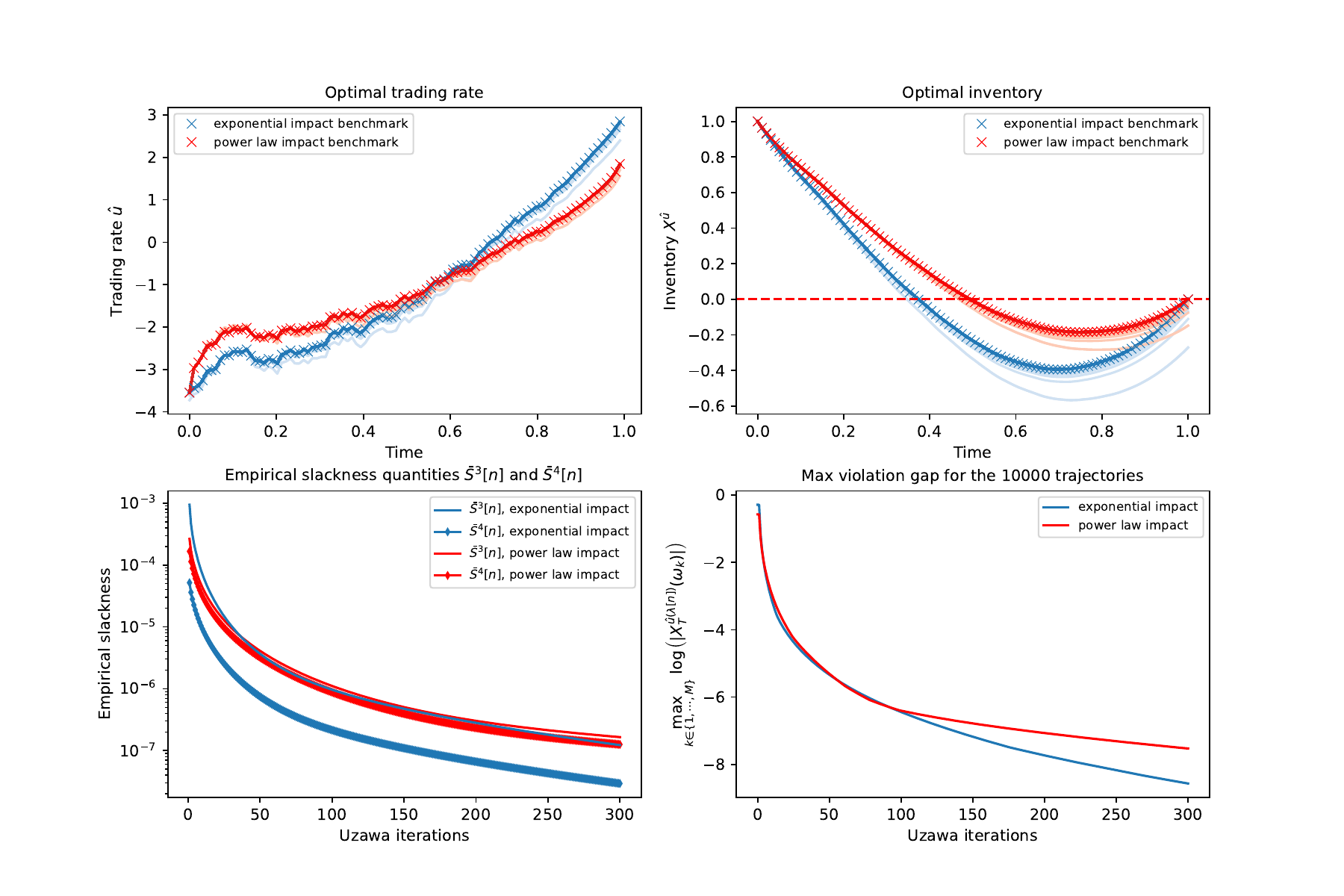}
    \vspace*{-0.5in}
    \caption{Illustration of path-wise convergence of the stochastic Uzawa Algorithm \ref{algo:stochatic_uzawa} for the optimal liquidation problem formulation \eqref{eq:optimal_liquidation} constructing $M = 1e4$ sample trajectories of the optimal control \eqref{eq:explicit_optimal_control}. Upper-left: one randomly selected optimal trading rate sample trajectory respectively for the exponential and the power law transient impact decays (random element being fixed); upper-right: their associated sample running inventories: we see both sample paths converging through Uzawa iterations (in plain line, with the color intensity increasing through the iterations) to the quadratic penalization benchmark \eqref{eq:benchmark_soft_constraint} (corresponding to the crosses) when $\nu = 1e6$ in \eqref{eq:benchmark_soft_constraint}. Lower-left: empirical slackness variables from \eqref{eq:discrete_time_complementary_slackness}; lower-right: maximum final inventory violation gap over all $1e4$ trajectories. The Uzawa algorithm is run with $300$ iterations, $\delta = 3$ and $\beta = 0.6$ in \eqref{eq:adaptive_learning_step} for the learning steps. Transient impact kernel's parameters for the exponential decay \eqref{eq:def_exponential_kernel} are set as: $c=5$ and $\rho=1$, and for the power law decay \eqref{eq:def_fractional_kernel}: $c=2$ and $\alpha = 0.6$. Signal parameters from \eqref{eq:drift_signal_specification}-\eqref{eq:seasonal_level_specification} are set as: $\theta=-20$, $w=0$, $\phi=\frac{\pi}{2}$, $\kappa = 1$, $\xi = 4$, $I_0 = -2$.}
    \label{F:convergence_plot}
    \end{center}
    \end{figure}

    \vspace{-0.65cm}

    \begin{remark}[Slater not restrictive for Stochastic Uzawa] \label{R:Slater_and_numerics}
        The sufficient condition of optimality as stated in Theorem \ref{T:generalized_kkt_theorem_sufficient} holds true for both equality and inequality constraints indifferently. On the other hand, the Slater condition $int \left( \mathcal{U} \right) \neq \emptyset$ used in Theorem \ref{T:generalized_kkt_theorem_necessary} excludes the case of equality constraints. However, it is worth noting that such assumption is not required to obtain the theoretical existence of Lagrange multipliers but only to ensure the separating hyperplane is non-vertical (see section \ref{ss:existence_lagrange_linear_form}). For example, the full liquidation requirement from \eqref{eq:assanity} indeed assumes an equality constraint at maturity $T$ on the inventory such that the Slater condition is not satisfied. In our numerical illustrations, we could have considered instead a relaxed formulation by setting the constraining functions as
        \begin{equation} \label{eq:assanity}
             -a_{t}^{1} = a_{t}^{2} = \bar M,  \quad t \in [0,T], \quad  -a_{t}^{3} =  a_{t}^{4} = \bar M, \quad t \in [0,T),\quad 
             - a_{T}^{3} = a_{T}^{4} = \epsilon, \quad \epsilon > 0,
        \end{equation}
        so that the Slater condition is satisfied and Theorem \ref{T:generalized_kkt_theorem_necessary} can be applied. However this is in practice unnecessary and in our numerical examples, no additional numerical stability was gained whatsoever by such $\epsilon$-relaxation formulation, so we set $\epsilon = 0$.
    \end{remark}
    
    \newpage
    \item[(ii)] \textbf{No buying with terminal liquidation.} Figure \ref{F:liquidation_with_buy_signal_example} compares $10$ sample trajectories (out of $1e4$) constructed by the stochastic Uzawa Algorithm \ref{algo:stochatic_uzawa} when the agent has to liquidate her inventory at maturity, in the presence of a positive `buy' signal, respectively with (in blue) and without (in red) a `no-buying' constraint. In this case, we specify the respective constraining functions from \eqref{eq:optimal_liquidation_admissible_set} as in \eqref{eq:assanity} for the agent who only aims at liquidating, and as
    \begin{align}\label{eq:no_buy_constraint}
         a_{t}^{1} = - \bar M, \; a_{t}^{2} = 0, \quad t \in [0,T], \quad  -a_{t}^{3} =  a_{t}^{4} =\bar M, \quad t \in [0,T),\quad 
         a_{T}^{3} = a_{T}^{4} = 0
    \end{align}
    for the agent who faces additionally a `no-buying' constraint (in blue). We also assume both agents face the same exponentially decaying transient market impact. Note that both traders' upper bound's Lagrange multiplier $\lambda^{4}$ got activated at terminal date $T$ in order to ensure the full liquidation constraint. The trader without the `no-buying' constraint (in red) capitalizes on the positive `buy' signal (in green) by initially increasing her inventory to profit from the expected price rise. She then liquidates rapidly to achieve a zero inventory by the terminal date.  
    
    In contrast, the trader facing the `no-buying' constraint (in blue) is not allowed to buy the asset. Consequently, she waits until about halfway through the trading horizon before beginning to liquidate, doing so at a  slower pace than the unconstrained trader. Notice that the Lagrange multiplier $\lambda^{2}$ is activated on the first half of the liquidation period for the constrained trader reflecting the `no-buying' constraint is indeed binding.

    \begin{figure}[H]
    \begin{center}
    \hspace*{-0.3in}
    \includegraphics[width=7 in,angle=0]{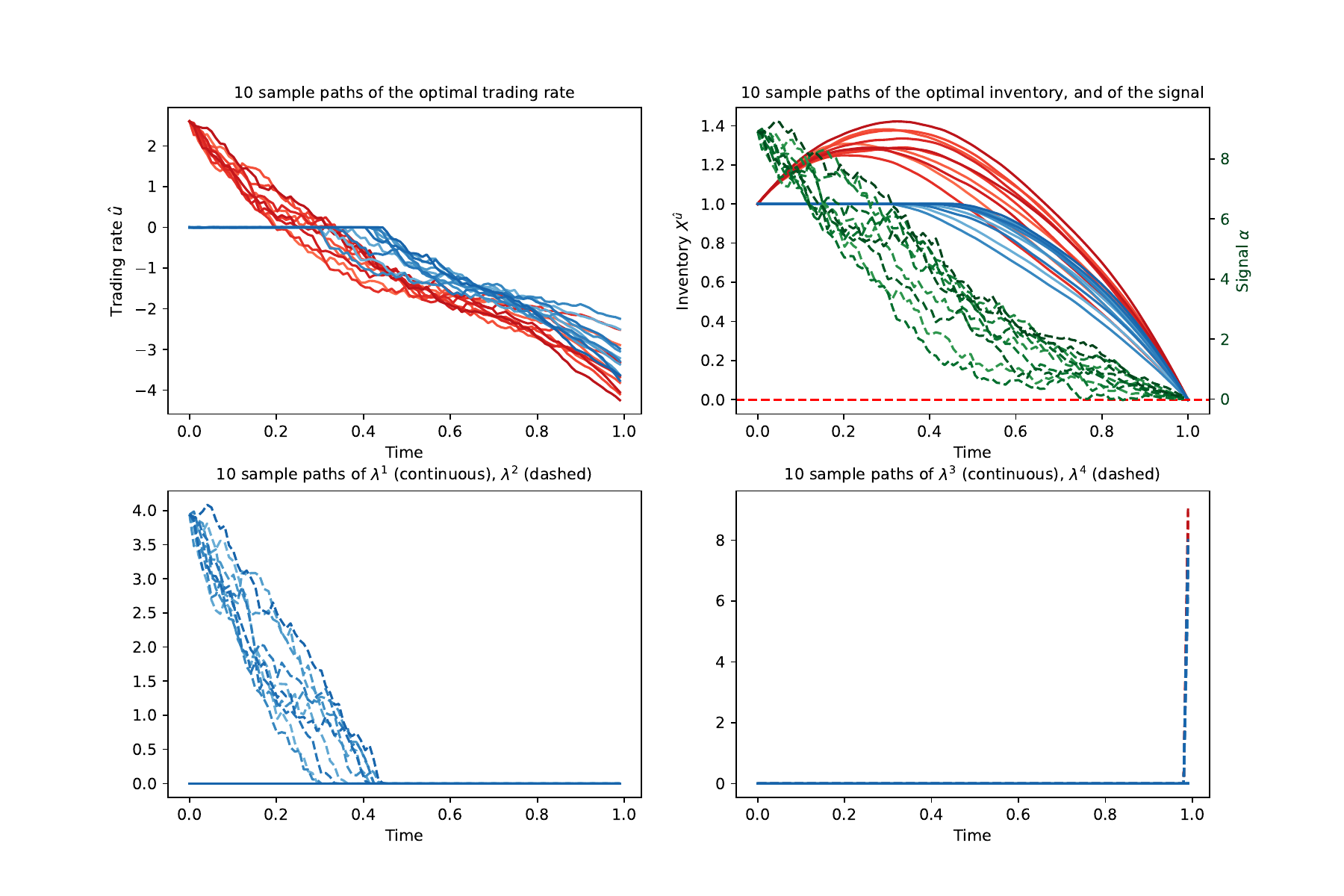}
    \vspace*{-0.5in}
    \caption{Optimal liquidation only (in red) including additionally the `no-buy' constraint (in blue) in presence of a `buy' $\alpha$-signal (in green, upper-right). The respective Lagrange multipliers are displayed in the bottom-plots using the same color code. Transient impact kernel's parameters for the exponential decay \eqref{eq:def_exponential_kernel} are fixed to $c=5$ and $\rho=1$. Signal parameters from \eqref{eq:drift_signal_specification}-\eqref{eq:seasonal_level_specification} are set as: $\theta=-5$, $w=0$, $\phi=\frac{\pi}{2}$, $\kappa = 1$, $\xi = 4$, $I_0 = 17$. Uzawa Algorithm \ref{algo:stochatic_uzawa} is run on $M=1e4$ trajectories, with $50$ iterations with learning steps specified in \eqref{eq:adaptive_learning_step} with $\delta = 1$, $\beta = 1e-4$ for the full liquidation controls (in red), and $50$ iterations with $\delta = 1$ and $\beta = 1e-4$ in \eqref{eq:adaptive_learning_step} for the full liquidation and `no-buy' constrained controls (in blue).}
    \label{F:liquidation_with_buy_signal_example}
    \end{center}
    \end{figure}

    \item[(iii)] \textbf{No shorting with terminal liquidation.} Figure \ref{F:liquidation_with_short_signal_example} compares $10$ sample trajectories (out of $1e4$) constructed by the stochastic Uzawa Algorithm \ref{algo:stochatic_uzawa} when the agent has to liquidate her inventory at maturity, in the presence of a `short signal', respectively with (in blue) and without (in red) a `no-shorting' constraint. In this case, the respective constraining functions from \eqref{eq:optimal_liquidation_admissible_set} are again specified as in \eqref{eq:assanity} for the agent who only aims at liquidating, and as
    \begin{align}\label{eq:no_short_constraint}
         - a_{t}^{1} = a_{t}^{2} = \bar M, \quad t \in [0,T], \quad a_{t}^{3} = 0, \; a_{t}^{4} = \bar M, \quad t \in [0,T),\quad 
         a_{T}^{3} = a_{T}^{4} = 0
    \end{align}
    for the agent who faces additionally a `no-short' constraint (in blue). We also assume both agents face the same exponentially decaying transient market impact. The trader without constraints liquidates rapidly in response to the `sell' signal (in green), aiming to complete a round-trip by shorting the asset before maturity. This strategy allows her to profit from the expected price decrease, which can significantly offset the price impact. She benefits further from the signal by holding the short position until the last moment, buying back the asset just before the terminal date to ensure a zero inventory. Observe that her Lagrange multiplier $\lambda^{3}$ only get activated at the terminal date $T$ to enforce the full liquidation constraint. In contrast, the trader facing the 'no-shorting' constraint cannot short the asset to profit from a round-trip strategy. Instead, she liquidates her inventory at a slower pace to minimize market impact. Once her inventory reaches zero, she ceases trading, as indicated by the early activation of her lower bound's Lagrange multiplier $\lambda^{3}$ before the terminal date $T$, ensuring she stops trading thereafter.

    \begin{figure}[H]
    \vspace*{-0.1in}
    \begin{center}
    \hspace*{-0.3in}
    \includegraphics[width=7 in,angle=0]{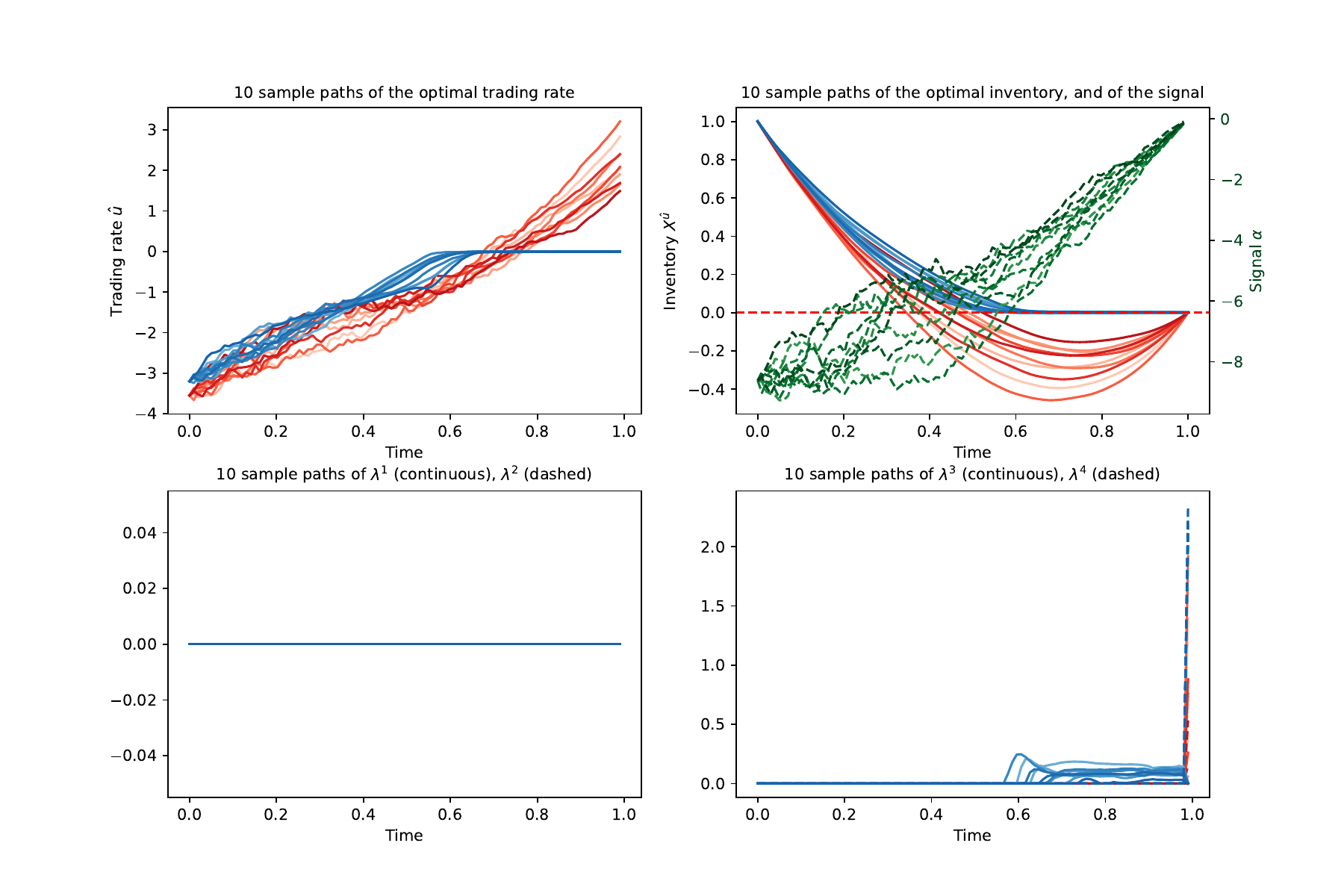}
    \vspace*{-0.5in}
    \caption{Optimal liquidation only (in red) and including additionally the no-short constraint (in blue) in presence of a "short" $\alpha$-signal (in green, upper-right). The respective Lagrange multipliers are displayed in the bottom-plots using the same color code. Transient impact kernel's parameters for the exponential decay \eqref{eq:def_exponential_kernel} are set as: $c=5$ and $\rho=1$. Signal parameters from \eqref{eq:drift_signal_specification}-\eqref{eq:seasonal_level_specification} are set as: $\theta=-20$, $w=0$, $\phi=\frac{\pi}{2}$, $\kappa = 1$, $\xi = 4$, $I_0 = -2$. Uzawa Algorithm \ref{algo:stochatic_uzawa} is run with $50$ iterations with learning steps specified in \eqref{eq:adaptive_learning_step} with $\delta = 3$, $\beta = 0.6$ for the full liquidation controls, and with $1e4$ iterations with $\delta = 0.1$ and $\beta = 5e-4$ in \eqref{eq:adaptive_learning_step} for the full liquidation and no-short constrained controls.}
    \label{F:liquidation_with_short_signal_example}
    \end{center}
    \end{figure}

    \newpage
    \item[(iv)] \textbf{Stop trading if the exogenous price $S$ drops below an a priori given reference price $S_{r} < S_{0}$.} Define the first-hitting process when $S$ drops below $S_{r}$ such that
    \begin{equation} \label{eq:stopping_time_drop_below_S_r}
        \tau^{r} := \inf_{t \in [0,T]} \left\{ S_{t} < S^{r} \right\} \wedge T,
    \end{equation}
    and we specify the constraining functions from \eqref{eq:optimal_liquidation_admissible_set} as
    \begin{equation} \label{eq:stop_trading_constraint}
         - a_{t}^{1} = a_{t}^{2} = \bar M \mathbb{1}_{[0, \tau_{r}]}(t), \quad t \in [0,T], \quad - a_{t}^{3} = a_{t}^{4} = \bar M', \quad t \in [0,T), \quad 
         - a_{T}^{3} = a_{T}^{4} = \bar M' \mathbb{1}_{[0, T)}(\tau^{r}).
    \end{equation}
    Figure \ref{F:stochastic_target} illustrates three sample trajectories of the stochastic optimal trading strategy, where we assume the agent faces a power law decaying transient market impact: 
    \begin{itemize}
        \item in red, the exogenous price $S$ never reaches the lower bound $S^{r}$ so that the agent liquidates her full inventory,
        
        \item in blue, the agent has to stop executing before reaching maturity as the `stop-trading' constraint got activated,

        \item in green, the `stop-trading' constraint got activated too, but the trader succeeds to unwind her position before being actually constrained to stop the execution.
    \end{itemize}

    Also, define the empirical average optimal trading rate by
    \begin{equation} \label{eq:average_trading_rate}
        \bar{u} := \frac{1}{M} \sum_{m = 1}^{M} \hat{u}(\omega_{m}),
    \end{equation}
    we observe that the average trading strategy (in magenta) indeed liquidates the inventory at maturity. 

    \begin{figure}[H]
    \vspace*{-0.15in}
    \begin{center}
    \hspace*{-0.3in}
    \includegraphics[width=7.2 in,angle=0]{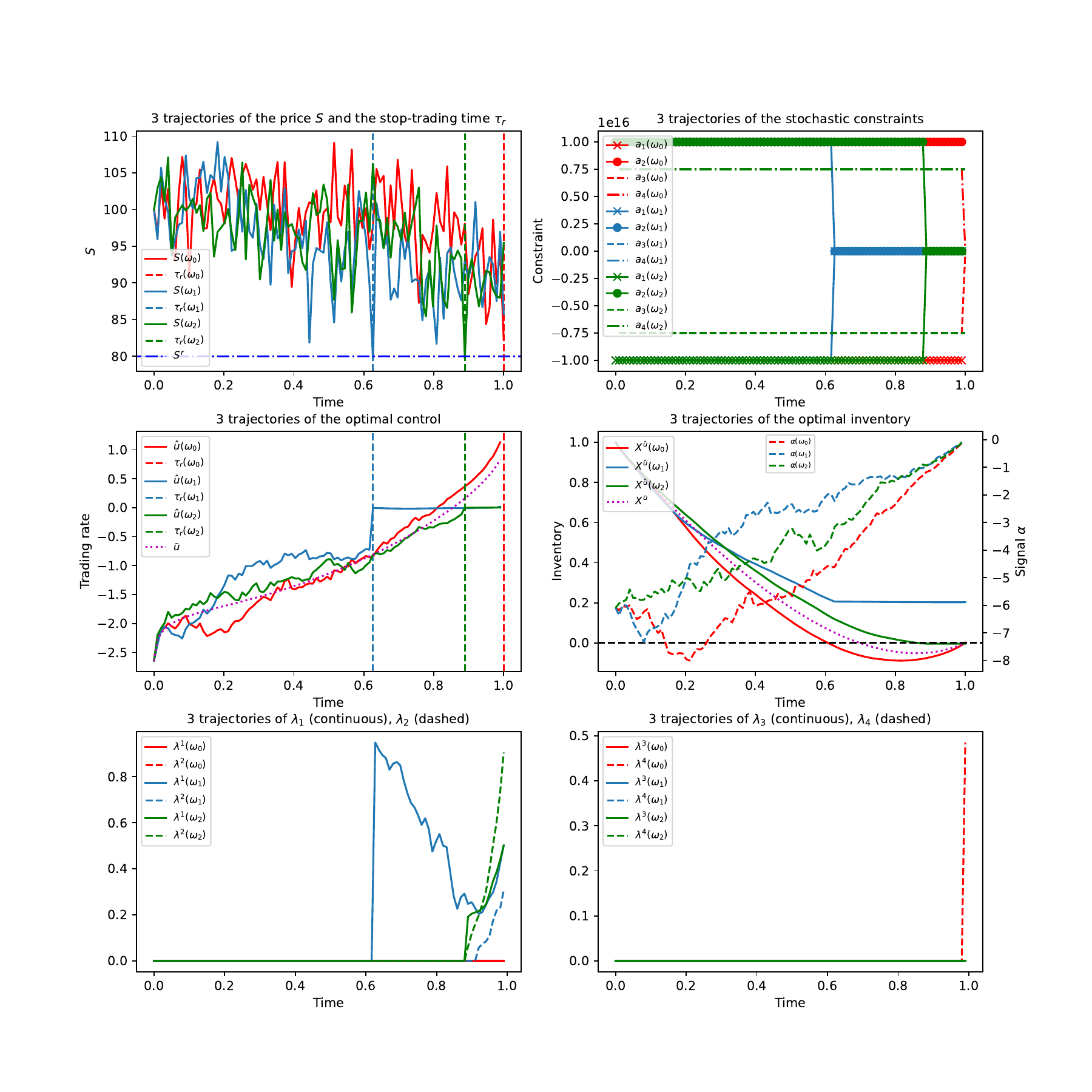}
    \vspace*{-0.5in}
    \caption{Optimal execution with `stop-trading' constraint if the exogenous price drops below a pre-specified level $S^{r} = 80$. Upper-left: $3$ sample trajectories of the exogenous price process $S$ from \eqref{eq:drift_like_signal_spec} with volatility $\sigma = 2$ and of the associated stopping time $\tau^{r}$ from \eqref{eq:stopping_time_drop_below_S_r}. Upper-right: $3$ sample trajectories of the corresponding constraining functions, where in both blue and green, the agent has to stop trading before reaching maturity; notice that for the sake of clarity, we specified $\bar M = 1e16$ and $\bar M' = 7.5e15$ in \eqref{eq:stop_trading_constraint}. Middle-left: the corresponding $3$ sample trading strategies. Middle-right: the corresponding $3$ sample running inventories and signals. Bottom left: the corresponding $3$ sample Lagrange multipliers for rate constraints. Bottom right: the corresponding $3$ sample Lagrange multipliers for inventory constraints. Transient impact kernel's parameters for the power-law decay \eqref{eq:def_fractional_kernel} are set to $c=2$ and $\alpha=0.6$. Signal parameters from \eqref{eq:drift_signal_specification}-\eqref{eq:seasonal_level_specification} are set as: $\theta=-20$, $w=0$, $\phi=\frac{\pi}{2}$, $\kappa = 1$, $\xi = 4$, $I_0 = -2$. For the selected model parameters, $194$ trajectories out of $1e4$ (i.e. $2\%$) get the `stop-trading' constraint activated. Uzawa Algorithm \ref{algo:stochatic_uzawa} is run on $1e4$ trajectories, with $5e3$ iterations with learning steps specified in \eqref{eq:adaptive_learning_step} with $\delta = 0.75$ and $\beta = 0.7$.}
    \label{F:stochastic_target}
    \end{center}
    \end{figure}

\end{enumerate}

\vspace{-0.3in}

\subsection{Energy storage with seasonality} \label{ss:energy_storage_example}

\begin{figure}[H]
\begin{center}
\hspace*{-0.3in}
\includegraphics[width=7 in,angle=0]{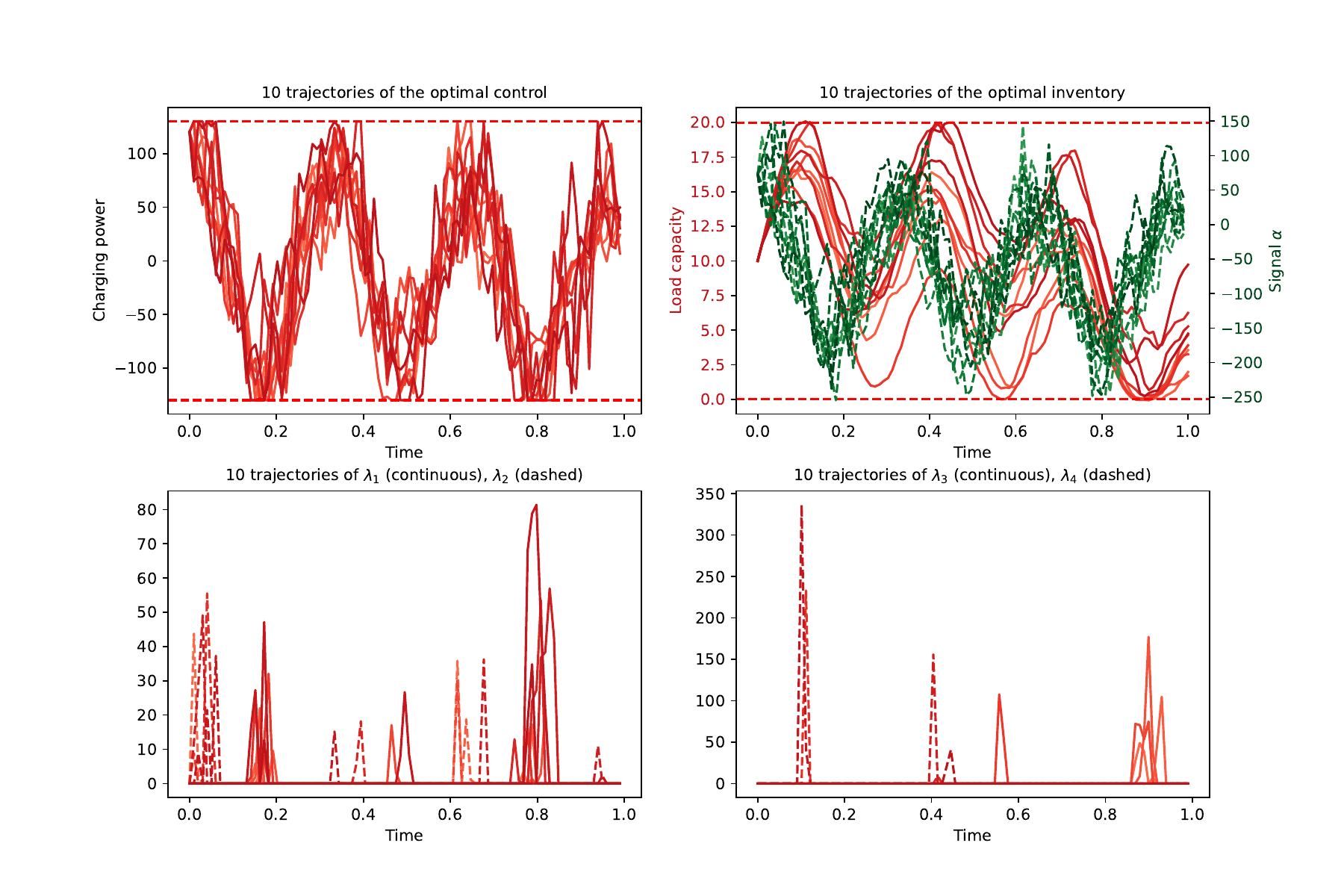}
\vspace*{-0.5in}
\caption{Optimal charging and associated load capacity trajectories (in red, upper-plots) in presence of a seasonal `buy' and `sell' $\alpha$-signal (in green, upper-right), with $u^{max}=-u^{min}=120$ (power unit) and $X^{min}=0$, $X^{max}=20$ (energy unit). The respective Lagrange multipliers are displayed in the bottom-plots. Transient impact parameters for the exponential kernel from \eqref{eq:def_exponential_kernel} are set to $c=5$ and $\rho=1$. Signal parameters from \eqref{eq:drift_signal_specification}-\eqref{eq:seasonal_level_specification} are set as: $\theta=1e5$, $w=20$, $\phi=0$, $\kappa = 5e1$, $\xi = 2e4$, $I_0 = 2e3$. Uzawa Algorithm \ref{algo:stochatic_uzawa} is run on $M=1e3$ trajectories, with $1e4$ iterations with learning steps specified as $\delta = 0.2$ and $\beta = 1e-3$ in \eqref{eq:adaptive_learning_step}.}
\label{F:battery_storage_example}
\end{center}
\end{figure}

Finally, we solve an optimal storage problem with constant constraints on both the charging power and the load capacity of the battery, formulated as
\begin{align} \label{eq:continuous_battery}
    \sup_{u \in \mathcal{U}} & \; \mathbb{E} \left[ \int_{0}^{T} \left( \alpha_{t} - \frac{1}{2} u_{t} \right) u_{t} \d t \right], \\
    & \mathcal{U} := \left\{ u \in \mathcal{L}^{\infty}, \; \; -u^{max} \leq u_{t} \leq u^{max} \; \text{and} \; 0 \leq X_{t}^{u} \leq X^{max} \right\},
\end{align}
where $u^{max}, X^{max} \in \mathbb R_{+}^{*}$. Figure~\ref{F:battery_storage_example} shows sample trajectories of the resulting optimal control obtained by the stochastic Uzawa Algorithm \ref{algo:stochatic_uzawa} when facing both seasonality and uncertainty on the price signal i.e.~with drift-like signal given by \eqref{eq:drift_signal_specification}--\eqref{eq:seasonal_level_specification} with a non-zero pulsation $w \neq 0$. Notice the respective Lagrange multipliers' trajectories become periodically positive when the constraints are indeed binding.

\newpage
\section{Existence and uniqueness of an optimal control}\label{S:existence}
In preparation for the proofs of our main theorems in Section~\ref{S:problem_formulation}, we start by proving the existence and uniqueness of an optimal control to the problem \eqref{eq:problem_formulation}.

\begin{theorem} \label{T:existenceuniqueness} Assume $\mathcal{U}$ is non-empty and let $K \in \mathcal K$ be an admissible kernel. There exists a unique optimal admissible control $\hat u \in \mathcal U$ such that $\mathcal J(\hat u)=\inf_{u\in\mathcal U} \mathcal J(u)$. 
\end{theorem}

The rest of the section is dedicated to the proof of Theorem~\ref{T:existenceuniqueness} which follows from general results on convex optimization in infinite dimensional spaces that we start by recalling. 

\subsection{A reminder on  convex control problems in $\mathcal{L}^{\infty}$} \label{ss:existence_uniqueness}

\noindent Existence results for infinite-dimensional convex control problems of the form
\begin{equation} \label{eq:generic_control_problem}
    \inf_{\mathcal{A}} \mathcal{T}
\end{equation}
where the set $\mathcal{A} \subset \mathcal{L}^{\infty}$ and $\mathcal{T}:\mathcal{L}^{\infty} \to \mathbb R$ are both convex\footnote{We also recall that the set $\mathcal{A}$ and the functional $\mathcal{T}$ are convex if $\theta x + (1-\theta) y \in \mathcal{U} \text{ and } \mathcal{T} \left( \theta x + (1-\theta) y \right) \leq \theta \mathcal{T}(x) + (1-\theta) \mathcal{T}(y), \quad x,y \in \mathcal{A}, \quad \theta \in (0,1)$.}, and satisfy
\begin{equation} \label{eq:proper_def}
    \mathcal{A} \cap \left\{ u \in \mathcal{L}^{\infty} | \mathcal{T}(u) < + \infty \right\} \neq \emptyset,
\end{equation}
typically rely the weak-$\star$ lower semi-continuity of the cost functional $\mathcal{T}$ such that
\begin{equation} \label{eq:lower_semi_continuity_def}
    \liminf_{f \to f_{0}} \mathcal{T}(f) \geq \mathcal{T}(f_{0}), \quad f_{0} \in \mathcal{L}^{\infty},
\end{equation}
where the limit is taken with respect to the weak-$\star$ topology $\sigma \left( \mathcal{L}^{\infty}, \mathcal{L}^{1} \right)$. Recall that the weak-$\star$ topology  is the coarsest topology that makes the family of linear forms indexed by $f \in \mathcal{L}^{\infty}$
\begin{align}
    \varphi_{f} : \; & \mathcal{L}^{1} \to \mathbb{R} \\
    & g \mapsto \varphi_{f}(g) := \int_{\Omega \times [0,T]} gf \left( \d \mathbb{P} \otimes \d t \right)
\end{align}
continuous, see \cite[Section 3.2]{brezis2011functional}. Three standard and useful results related to this topology include:
\begin{enumerate}
    \item[(i)] every bounded sequence valued in $\mathcal{L}^{\infty}$ admits a converging sub-sequence for the weak-$\star$ topology $\sigma \left( \mathcal{L}^{\infty}, \mathcal{L}^{1} \right)$ \cite[Corollary 3.30]{brezis2011functional} due to the fact that it is the dual of $\mathcal{L}^{1}$, which is indeed a separable Banach space given that $\Omega \times [0,T]$ is separable \cite[Theorems 4.14 and 4.13]{brezis2011functional},
    
    \item[(ii)] any convex and closed subset of $\mathcal{L}^{\infty}$ is also weakly-$\star$-closed \cite[Theorem 3.7]{brezis2011functional},

    \item[(iii)] Banach–Alaoglu–Bourbaki's Theorem  \cite[Theorem 3.16]{brezis2011functional} ensures that bounded weak-$\star$-closed sets are weak-$\star$ compact.
\end{enumerate}

Finally, for $\beta > 0$, we say that
$\mathcal{T}$ is $\beta$-convex if
\begin{equation} \label{eq:beta_convexity}
     \mathcal{T} (\theta u + (1-\theta) v) \leq \theta \mathcal{T}(u) + (1-\theta) \mathcal{T}(v) - \frac{\beta \theta (1-\theta)}{2} ||u-v||_{2}^2, \quad (u,v) \in \left(\mathcal{L}^{\infty}\right)^2, \quad \theta \in (0,1).
\end{equation}
We are now ready to state an existence and uniqueness result.
\begin{theorem}[Existence (and uniqueness)] \label{T:existence_uniqueness}
    Assume \eqref{eq:proper_def}, and that $\mathcal{A}$ is a closed, bounded\footnote{A weaker condition of coercivity of the functional of the form $\lim_{\| u \|_{\infty} \to \infty} \mathcal{T}(u) = \infty$ can replace the assumption of $\mathcal{A}$ being bounded.} and convex non-empty subset of $\mathcal{L}^{\infty}$, that $\mathcal{T}$ is a real-valued, convex and lower-semi-continuous functional on $\mathcal{L}^{\infty}$. Then, the general convex program
    \begin{equation}
        \inf_{u \in \mathcal{A}} \mathcal{T}(u)
    \end{equation}
    admits a solution. Furthermore, if  $\mathcal{T}$ is $\beta$-convex for some $\beta>0$, then the minimizer of $\mathcal{T}$ on $\mathcal{A}$ is unique, denoted by $\hat{u}$, and satisfies
    \begin{equation} \label{eq:error_estimate}
        || u - \hat{u} ||_{2}^{2} \leq \frac{4}{\beta} (\mathcal{T}(u) - \mathcal{T}(\hat{u})), \quad u \in \mathcal{A}.
    \end{equation}
\end{theorem}

\begin{proof}
    First let us prove the existence. Consider a minimizing sequence $(u^{n}) \in \mathcal{A}^{\mathbb{N}}$ such that
    \begin{equation}
        \lim_{n \to \infty} \mathcal{T}(u^{n}) = \inf _{u \in \mathcal{A}} \mathcal{T}(u),
    \end{equation}
    and which is a bounded since $\mathcal{A}$ is bounded. Consequently, there is a sub-sequence $(u^{n_{k}})$ that converges for the weak-$\star$ topology $\sigma \left( \mathcal{L}^{\infty}, \mathcal{L}^{1} \right)$ to a limit $\hat{u}$, cf.point (i) above, and since $\mathcal{A}$ is convex and closed, it is also weak-$\star$-closed and we get $\hat{u} \in \mathcal{A}$, cf.points (ii)-(iii). Finally, the weak-$\star$ lower-semi-continuity \eqref{eq:lower_semi_continuity_def} of $\mathcal{T}$ yields
    \begin{equation}
        \mathcal{T}(\hat{u}) \leq \liminf_{k \to \infty} \mathcal{T}(u^{n_{k}}) = \inf_{v \in \mathcal{A}} \mathcal{T}(v),
    \end{equation}
    where the equality is given by definition of the minimizing sequence, which yields the existence of a solution.   
    Uniqueness follows by strict convexity: assume there are two solutions $\hat{u}^{1}$ and $\hat{u}^{2}$ such that
    \begin{equation}
        \inf_{v \in \mathcal{A}} \mathcal{T}(v) \leq \mathcal{T} \left( \frac{1}{2}\hat{u}^{1} + \frac{1}{2}\hat{u}^{2} \right) < \frac{1}{2}\mathcal{T}(\hat{u}^{1}) + \frac{1}{2}\mathcal{T}(\hat{u}^{2}) = \inf_{v \in \mathcal{A}} \mathcal{T}(v),
    \end{equation}
    which is indeed a contradiction. Finally, take $\theta = \frac{1}{2}$ in the definition \eqref{eq:beta_convexity} of $\beta$-convexity to get
    \begin{equation} \label{eq:interm_error_bound}
        \frac{\beta}{8} \| u - \hat{u} \|_{2}^{2} \leq \frac{1}{2} \mathcal{T}(u) + \frac{1}{2} \mathcal{T}(\hat{u}) - \mathcal{T} \left( \frac{u+\hat{u}}{2} \right) \leq \frac{\mathcal{T}(u) - \mathcal{T}(\hat{u})}{2}
    \end{equation}
    where the second inequality holds because
    \begin{equation} \label{eq:interm_error_bound}
        \mathcal{T}(\hat{u}) \leq \mathcal{T} \left( \frac{u+\hat{u}}{2} \right),
    \end{equation}
    and re-arranging the terms yields \eqref{eq:error_estimate}. 
\end{proof}

\subsection{Proof of Theorem~\ref{T:existenceuniqueness}} \label{ss:proofexistenceuniqueness}

We start by reformulating the cost functional \eqref{eq:functional_formulation_storage} to make it easier to manipulate. 

\begin{lemma} \label{L:problem_reformulation_alpha_cost_trade_of}
The cost functional $\mathcal{J}$ in \eqref{eq:functional_formulation_storage} can be re-written as
    \begin{equation} \label{eq:functional_after_ipp}
        \mathcal{J}(u) = \mathbb E \left[ \int_{0}^T \left( - \alpha_{t} + \frac{\gamma}{2} u_t + \int_0^t K(t,s) u_s \d s\right) u_t \d t\right] - X_{0} \E \left[ S_T \right],
    \end{equation}
    with the signal process $\alpha$ given by \eqref{eq:alpha_signal_definition}.
\end{lemma}

\begin{proof}
    Using an integration by part yields
    \begin{equation} \label{eq:ipp_}
        \int_0^T S_t u_t \d t - X_{T}^u S_T = - \int_0^T \left(\int_0^t u_s \d s\right) \d S_t - X_{0} S_T,
    \end{equation}
    and taking the expectation of \eqref{eq:ipp_} while injecting the dynamics of $S$ given by \eqref{eq:exogenous_price_process} yields
    \begin{equation} \label{eq:expectation_ipp_}
        \E \left[ \int_0^T S_t u_t \d t + X_{T}^u S_T \right]= - \E \left[ \int_0^T \left( \int_0^t u_s \d s \right) \d P_t \right] - X_{0} \E \left[ S_T \right],
    \end{equation}
    where we used the fact that 
    \begin{equation}
        \E \left[ \int_0^T \left( \int_0^t u_s \d s \right) \d M_t \right] = 0,
    \end{equation}
    recall that $M$ is a centered martingale. Now, applying Fubini to \eqref{eq:expectation_ipp_},  using the tower property,   and injecting into \eqref{eq:functional_formulation_storage} leads to  the result. 
\end{proof}

\begin{lemma} \label{L:ppties_set_of_admissible_controls_and_functional}
Assume $\mathcal{U}$ is non-empty and let $K \in \mathcal K$ be an admissible kernel. The set $\mathcal U$ is a closed, bounded and convex subset of $\mathcal L^{\infty}$ and  the cost functional $\mathcal J$ is $1$-convex in the sense of \eqref{eq:beta_convexity} and continuous.
\end{lemma}

\begin{proof}
$\mathcal U$ is clearly closed, bounded and convex in $\mathcal L^{\infty}$. To prove the strong convexity on $\mathcal J$ we start by decomposing the functional in \eqref{eq:functional_after_ipp} into
    \begin{equation}
        \mathcal{J} := \mathcal{J}_{1} + \mathcal{J}_{2},
    \end{equation}
    with
    \begin{equation}
        \mathcal{J}_{1}(u) := \E \left[ \int_{0}^T \left(- \alpha_t + \frac{1}{2} u_{t} \right) u_{t} \d t \right] \quad \text{ and } \quad \mathcal{J}_{2}(u) := \E \left[ \int_{0}^T \int_0^t K(t,s) u_{s} \d s u_{t} \d t \right],
    \end{equation}
    and proving that $\mathcal{J}_{1}$ is $1$-convex and $\mathcal{J}_{2}$ is convex. Let $\theta \in (0,1)$, $u,v \in \mathcal{L}^{\infty}$ and observe that by completing the square
    \begin{equation}
        \mathcal{J}_{1} \left( \theta u + (1-\theta) v \right) = \theta \mathcal{J}_{1} \left( u \right) + (1-\theta) \mathcal{J}_{1} \left( v \right) - \frac{\theta (1 - \theta)}{2} \left| \left| u-v \right| \right|_{2}^{2} + \frac{\theta (\theta - 1)}{2} \E \left[ \int_{0}^{T} \left( u_{t}^{2} + v_{t}^{2} \right) \d t \right],
    \end{equation}
    which readily yields that $\mathcal{J}_{1}$ is indeed $1$-convex in the sense of \eqref{eq:beta_convexity}, since the last term is non-negative. Finally, the convexity of $\mathcal{J}_{2}$ follows from the same reasoning as in the proof of \cite[Lemma 5.1]{abijaberOptPort2024}, using the non-negative definite property of the kernel $K$ from \eqref{eq:pos_def} and it is clear that $J$ is continuous on $\mathcal L^{\infty}$: it is even Fréchet-differentiable, cf. Lemma \ref{L:properties_functional_admin_controls}.
\end{proof}

We can now prove Theorem~\ref{T:existenceuniqueness}.
\begin{proof}[Proof of Theorem~\ref{T:existenceuniqueness}]
     Using Lemma \ref{L:ppties_set_of_admissible_controls_and_functional}, we know in particular that the cost functional $\mathcal{J}$ is convex and continuous implying that it is weak-$\star$ lower-semi-continuous, while all the other assumptions of Theorem~\ref{T:existence_uniqueness} are readily satisfied, yielding existence and uniqueness of the optimal control $\hat u \in \mathcal U$.
\end{proof}

\begin{remark}
    Notice that the functional $\mathcal{J}$ without instantaneous market impact, i.e.~$\gamma=0$, is convex and the control problem still admits a solution, which may not be unique. Clearly, slippage costs in the model have a regularizing effect.
\end{remark}

\section{Proof of Theorems \ref{T:generalized_kkt_theorem_necessary} and \ref{T:generalized_kkt_theorem_sufficient}} \label{ss:proof_theorem_necessary_sufficient_condition}

In all of this section, we assume $\mathcal{U}$ is non-empty and let $K \in \mathcal K$ be an admissible kernel. We let  $\hat{u}$ be the unique minimizer of \eqref{eq:stochasticConvexProgramFormulation} obtained from Theorem~\ref{T:existenceuniqueness}. We denote by $\mathcal{Y} := \left( \mathcal{L}^{\infty} \right)^{4}$ which is a Banach space.

\subsection{$\mathcal J$ and $\mathcal G$ are Fréchet-differentiable}

\begin{lemma}
\label{L:properties_functional_admin_controls}
Both $\mathcal{J}$ and $\mathcal{G}$ respectively from \eqref{eq:stochasticConvexProgramFormulation} and \eqref{eq:def_constraining_functional} are Fréchet-differentiable at any  point $u \in \mathcal{L}^{\infty}$ with respective differentials $\nabla \mathcal{J}$ and $\nabla \mathcal{G}$ given for any $h \in \mathcal{L}^{\infty}$ by
    \begin{align} \label{eq:nabla_J}
        \nabla \mathcal{J}(u)(h) & := \E \left[ \int_0^T \left( - \alpha_t + u_t + \int_0^t K(t,s) u_s \d s + \int_t^T K(s,t) \E_t \left[ u_s \right] \d s \right) h_t \d t \right] \in \mathbb{R}, \\ \label{eq:nabla_G}
        \nabla \mathcal{G}(u)(h) & := \left( -h, h, - \int_0^. h_{s} \d s, \int_0^. h_{s} \d s \right) \in \mathcal{Y}.
    \end{align}
\end{lemma}

\begin{proof}
    \noindent Fix $u \in \mathcal{L}^{\infty}$ and $h \in \mathcal{L}^{\infty}$. We have by straightforward calculus
    \begin{equation}
        \mathcal{J}(u+h) = \mathcal{J}(u) + \nabla \mathcal{J}(u)(h) + \E \left[ \int_0^T \left( \frac{\gamma}{2} h_t + \int_0^t K(t,s) h_s \d s \right) h_t \d t \right],
    \end{equation}
    where $\nabla \mathcal{J}(u)(h)$ is given by \eqref{eq:nabla_J} and where we note that
    \begin{align}
        0 \leq \frac{\left| \E \left[ \int_0^T \left( \frac{\gamma}{2} h_t + \int_0^t K(t,s) h_s \d s \right) h_t \d t \right] \right|}{||h||_{\infty}} \leq \left( \frac{\gamma}{2} + \int_{0}^{T} \int_{0}^{T} \left| K(t,s) \right| \d s \d t \right) ||h||_{\infty} \to 0, \text{ as } ||h||_{\infty} \to 0,
    \end{align}
    where we used the triangular inequality, the positivity of the integral and the fact $K$ satisfies \eqref{eq:square_integrable_kernel_assumption} to conclude on the limit. The existence of the Fréchet differential of $\mathcal{G}$ is immediate noting that $\mathcal{G}(u+h) = \mathcal{G}(u) + \nabla \mathcal{G}(u)(h)$, with $\nabla \mathcal{G}(u)(h)$ given by \eqref{eq:nabla_G}.
\end{proof}

\subsection{The Lagrange stationary equation} \label{ss:existence_lagrange_linear_form}

We prove in this section the existence of a Lagrange continuous linear form associated to the constraining functional $\mathcal{G}$ defined in \eqref{eq:def_constraining_functional} and obtain an associated Lagrange stationary equation in functional form as a necessary condition for optimality. We will simplify and adapt the proofs of \cite[Theorem 5.3 and Corollary 5.4]{jahn1994introduction} to fit our problem and refer the reader to this reference for broader context where the minimization of a general functional (not necessarily convex) under both equality and inequality constraints is considered.

\begin{theorem}[Existence of a Lagrange separating linear form] \label{T:existence_lagrange_separating_functional}
    There exists a continuous linear functional $l \in \mathcal{Y}^{\star} \backslash \left\{ 0_{\mathcal{Y}^{\star}} \right\}$ such that
    \begin{equation} \label{eq:stationary_lagrange}
        \nabla \mathcal{J}(\hat{u}) + l \circ \nabla \mathcal{G} (\hat{u}) = 0_{(\mathcal{L}^{\infty})^{\star}},
    \end{equation}
    and
    \begin{equation} \label{eq:complementary_slackness_functional}
        l(\mathcal{G}(\hat{u})) = 0
        ,
    \end{equation}
    and
    \begin{equation} \label{eq:non_negative_separating_linear_form}
        l\left( y \right) \geq 0 \quad y \in \mathcal{C},
    \end{equation}
    i.e.~$l$ is an element of the dual cone $\mathcal{C}^{*}$.
\end{theorem}

The rest of the subsection is dedicated to the proof of Theorem~\ref{T:existence_lagrange_separating_functional}. First, let us recall the notions of \textit{cone} and \textit{generated cone} which will be useful in the following.

\begin{definition}[Cone and generated cone]
    Let $Q$ and $S$ be nonempty subsets of a real linear space. We say that $Q$ is a cone if
    \begin{equation}
        v \in Q \; \implies \; \lambda v \in Q, \quad \lambda \geq 0,
    \end{equation}
    and we define the cone generated by $S$ as
    \begin{equation}
        cone \left( S \right) := \left\{ \lambda s, \; \lambda \geq 0 \text{ and } s \in S \right\}.
    \end{equation}
\end{definition}

We now state a straightforward Lemma which will be useful later on.

\begin{lemma} \label{L:empty_intersection}
    There is no $u \in \mathcal{L}^{\infty}$ such that
    \begin{equation}
        \nabla \mathcal{J}(u-\hat{u}) < 0 \text{ and } \mathcal{G}(\hat{u}) + \nabla \mathcal{G}(\hat{u})(u-\hat{u}) \in -\text{int}(C).
    \end{equation}
\end{lemma}
\begin{proof}
    It is clear that if $u \in \mathcal{L}^{\infty}$ is such that
    \begin{equation}
        -\text{int}(C) \ni \mathcal{G}(\hat{u}) + \nabla \mathcal{G}(\hat{u})(u-\hat{u}) = \mathcal{G}(u),
    \end{equation}
    then it means $u \in \mathcal{U}$ and by the Euler-Lagrange optimality inequality, $u$ must also satisfy
    \begin{equation}
        \nabla \mathcal{J}(\hat{u})(u-\hat{u}) \geq 0, \quad u \in \mathcal{U},
    \end{equation}
    see for example \cite[Theorem 2.5.1]{allaire2023optimisation}.
\end{proof}

Let us recall the \textit{Eidelheit separation theorem} from functional analysis which actually will be the keystone in proving the existence of the Lagrange separating linear form.

\begin{theorem}[Eidelheit separation] \label{T:eidelheit_separation}
    Let $S$ and $T$ be nonempty convex subsets of a real topological linear space $X$ with $\text{int}(S) \neq \emptyset$. Then we have $\text{int}(S) \cap T = \emptyset$ if and only if there are a continuous linear functional $l \in X^{*} \backslash \{0_{X^{*}}\}$ and a real number $\xi$ with
    \begin{equation}
        l(s) \leq \xi \leq l(t), \quad s \in S, \quad t \in T,
    \end{equation}
    and
    \begin{equation}
        l(s) < \xi, \quad s \in \text{int}(S).
    \end{equation}
\end{theorem}

\begin{proof}
    See for example Theorem 3 from \cite[Section 5.12, p.133]{luenberger1997optimization}.
\end{proof}

We are now ready to state the existence theorem for a Lagrange linear form associated to the unique optimal solution $\hat{u}$ to our control problem.

\begin{proof}[Proof of Theorem~\ref{T:existence_lagrange_separating_functional}]
    Note that the set
    \begin{equation}
        Q := \left\{ (\nabla \mathcal{J}(\hat{u})(u-\hat{u})+\alpha, \mathcal{G}(\hat{u}) + \nabla \mathcal{G}(\hat{u})(u-\hat{u}) + c) \in \R \times \mathcal{Y} \; | \; u \in \mathcal{L}^{\infty}, \; \alpha > 0, \; c \in \text{int}(\mathcal{C}) \right\},
    \end{equation}
    is non-empty since $\mathcal{L}_{+}^{\infty}$ has a non-empty interior, and is convex by linearity of the Fréchet differentials $\nabla \mathcal{J}(\hat{u})$ and $\nabla \mathcal{G}(\hat{u})$.
    Then, an application of Lemma~\ref{L:empty_intersection} ensures
    \begin{equation}
        (0, 0_{\mathcal{Y}}) \notin Q, \quad i.e.~Q \cap (0, 0_{\mathcal{Y}}) = \emptyset.
    \end{equation}
    Note that $\R \times \mathcal{Y}$ is indeed a real Banach so in particular it is a real linear topological space (whose topology is induced by the metric associated to the product-norm) so that the \textit{Eidelheit separation} Theorem~\ref{T:eidelheit_separation} yields the existence of two real numbers $\mu, \; \xi$ and a linear functional $l \in \mathcal{Y}^{\star}$ with $(\mu, l) \neq (0,0_{\mathcal{Y}^{\star}})$ such that
    \begin{equation} \label{eq:ineq_intermediary_1}
        \mu \left( \nabla \mathcal{J}(\hat{u})(u-\hat{u})+\alpha \right) + l\left( \mathcal{G}(\hat{u}) + \nabla \mathcal{G}(\hat{u})(u-\hat{u}) + y \right) > \xi \geq 0, \quad (u, \alpha, y) \in \mathcal{L}^{\infty} \times \R_{+}^{*} \times \text{int}(\mathcal{C}).
    \end{equation}
    Noting that $\mathcal{C}$ is convex\footnote{See for example Theorem~4.3 from \cite{jahn1994introduction}.} and that every convex subset of a real normed space with nonempty interior is contained in the closure of the interior of this set, then we get from inequality \eqref{eq:ineq_intermediary_1}
    \begin{equation} \label{eq:ineq_intermediary_2}
        \mu \left( \nabla \mathcal{J}(\hat{u})(u-\hat{u})+\alpha \right) + l\left( \mathcal{G}(\hat{u}) + \nabla \mathcal{G}(\hat{u})(u-\hat{u}) + y \right) \geq \xi \geq 0, \quad (u, \alpha, y) \in \mathcal{L}^{\infty} \times \R_{+} \times \mathcal{C}.
    \end{equation}
    From the inequality \eqref{eq:ineq_intermediary_2}, we obtain for $u=\hat{u}$
    \begin{equation} \label{eq:ineq_intermediary_3}
        \mu \alpha + l\left( \mathcal{G}(\hat{u}) + y \right) \geq 0, \quad (\alpha, y) \in \R_{+} \times \mathcal{C}.
    \end{equation}
    With $\alpha=1$ and $y=-\mathcal{G}(\hat{u})$, we get $\mu \geq 0$. From the inequality \eqref{eq:ineq_intermediary_3}, it follows for $\alpha=0$
    \begin{equation} \label{eq:ineq_intermediary_4}
        l\left( \mathcal{G}(\hat{u}) \right) \geq - l\left( y \right), \quad y \in \mathcal{C}.
    \end{equation}
    Assume $l\left( y \right) < 0$ for some $y \in \mathcal{C}$, then with $\lambda y \in \mathcal{C}$, for some sufficiently large $\lambda > 0$, one gets a contradiction to the inequality \eqref{eq:ineq_intermediary_4}, yielding consequently the positivity of $l$ as in \eqref{eq:non_negative_separating_linear_form}. Moreover the inequality \eqref{eq:ineq_intermediary_4} implies that $l\left( \mathcal{G}(\hat{u}) \right) \geq 0$. Since $\hat{u}$ satisfies the inequality constraint, i.e.~$\mathcal{G}(\hat{u}) \in -\mathcal{C}$, we also get with inequality \eqref{eq:non_negative_separating_linear_form} that $l\left( \mathcal{G}(\hat{u}) \right) \leq 0$, hence we get $l\left( \mathcal{G}(\hat{u}) \right) = 0$, i.e.~the equation \eqref{eq:complementary_slackness_functional} is proven.

    Next, let us prove equation \eqref{eq:stationary_lagrange}. For $\alpha = 0$ and $y = -\mathcal{G}(\hat{u})$, we obtain from inequality \eqref{eq:ineq_intermediary_2}
    \begin{equation}
        \left( \mu \nabla \mathcal{J}(\hat{u}) + l \circ \nabla \mathcal{G}(\hat{u}) \right)(u-\hat{u}) \geq 0, \quad u \in \mathcal{L}^{\infty}.
    \end{equation}
    By linearity of both $\nabla \mathcal{J}(\hat{u})$ and $l$, we immediately deduce that
    \begin{equation} \label{eq:eq_intermediary}
        \mu \nabla \mathcal{J}(\hat{u}) + l \circ \nabla \mathcal{G}(\hat{u}) = 0_{(\mathcal{L}^{\infty})^{*}}.
    \end{equation}

    Given that $int(\mathcal{U}) \neq \emptyset$, then there is a control $\mathring{u}$ such that 
    $$
    \mathcal{G}(\mathring{u}) = \mathcal{G}(\hat{u}) + \nabla \mathcal{G}(\hat{u})(\mathring{u}-\hat{u}) \in - \text{int}(\mathcal{C}), 
    $$
    and an application of \cite[Theorem 5.6]{jahn1994introduction} yields the following Kurcyusz-Robinson-Zowe (KRZ) regularity condition
    \begin{equation} \label{eq:Kurcyusz_Robinsion_Zowe}
        \nabla \mathcal{G} (\hat{u}) \left( cone(\mathcal{L}^{\infty} - \{\hat{u}\}) \right) \; +\; cone(\mathcal{C}+\mathcal{G}(\hat{u})) = \mathcal{Y}.
    \end{equation}
    Let us now prove  that $\mu >0$. Given the KRZ regularity condition \eqref{eq:Kurcyusz_Robinsion_Zowe}, for an arbitrary $y \in \mathcal{Y}$, we can fix two non-negative real numbers $r_{1}, r_{2}$ as well as $u \in \mathcal{L}^{\infty}$ and $c \in \mathcal{C}$ such that
    \begin{equation}
        y = \nabla \mathcal{G}(\hat{u}) \left( r_{1}(u-\hat{u}) \right) + r_{2} \left( c + \mathcal{G}(\hat{u}) \right).
    \end{equation}
    Assume that $\mu = 0$. Then using the inequality \eqref{eq:eq_intermediary}, the equation \eqref{eq:complementary_slackness_functional} and the positivity of $l$ from \eqref{eq:non_negative_separating_linear_form}, we get by linearity
    \begin{equation}
        l(y) = r_{1} \left( l \circ \nabla \mathcal{G} (\hat{u}) \right) \left( u-\hat{u} \right) + r_{2} \left( l(c) + l\left( \mathcal{G}(\hat{u}) \right) \right) = r_{2} l(c) \geq 0,
    \end{equation}
    and, as a consequence of the linearity of $l$, we have $l=0_{\mathcal{Y}^{\star}}$ which contradicts $(\mu, l) \neq (0,0_{\mathcal{Y}^{\star}})$. Hence $\mu > 0$, and we can redefine $l := \frac{l}{\mu} \in \mathcal{C}^{*}$ such that we obtain indeed equality \eqref{eq:stationary_lagrange}.
\end{proof}

\vspace{1mm}

\subsection{The Lagrange multipliers and the Fredholm equation} \label{ss:lagrange_multipliers_and_fredholm}

In the next step, we will characterize the separating linear functional obtained in Theorem~\ref{T:existence_lagrange_separating_functional} in terms of  Lagrange multipliers lying in the $\textbf{\textit{ba}}$ space, and derive the stochastic Fredhlom equation satisfied by the optimal control $\hat{u}$.

\begin{lemma}\label{L:lagrangemultiplierandFredholm}
There exist Lagrange multipliers $\psi := \left( \psi^{1}, \psi^2,\psi^3,\psi^4 \right) \in \textbf{\textit{ba}}_{+}^{4}$ satisfying the equations \eqref{eq:aggregated_lagrange_multiplier} with $\varphi \in \mathcal{L}^{2}$ and \begin{equation} \label{eq:aggregated_lagrange_multiplier_singular}
     0 = \left( \psi^{1,o} - \psi^{2,o} + I^{*}\left(\psi^{3,o}\right) - I^{*}\left(\psi^{4,o}\right) \right) (\d t,\d \omega),
\end{equation}
as well as  the complementary slackness \eqref{eq:KKT_slackness_conditions} and \begin{equation} \label{eq:KKT_slackness_conditions_singular}
    \begin{cases}
       0 = \int_{[0,T] \times \Omega} \left(a_{t}^{1}(\omega) - \hat{u}_{t}(\omega)\right) \psi^{1,o}(\d t,\d \omega)  \\
       0 = \int_{[0,T] \times \Omega} \left(\hat{u}_{t}(\omega) - a_{t}^{2}(\omega)\right) \psi^{2,o}(\d t,\d \omega)  \\
       0 = \int_{[0,T] \times \Omega} \left( a^3_t (\omega) - X^{\hat{u}}_t(\omega) \right) \psi^{3,o}(\d t,\d \omega) \\
       0 = \int_{[0,T] \times \Omega} \left( X^{\hat{u}}_t(\omega) - a^4_t (\omega)  \right) \psi^{4,o}(\d t,\d \omega)
    \end{cases},
\end{equation}
and such that the stochastic Fredholm equation \eqref{eq:Fredholm_equation} is satisfied.
\end{lemma}

\begin{proof}
Fix $l \in \mathcal{Y}^{\star} \backslash \left\{ 0_{\mathcal{Y}^{\star}} \right\}$ as in Theorem~\ref{T:existence_lagrange_separating_functional}. Then
an application of \cite[Theorem 16, Chapter IV, Section 8]{dunford1988linear} yields the existence of four bounded and finitely additive set-valued functions $\left( \psi^{1}, \psi^2,\psi^3,\psi^4 \right) \in \textbf{\textit{ba}}^{4}$ such that
\begin{equation}
    l(x) = \sum_{i=1}^{4} \int_{\Omega \times [0,T]} x_{t}^{i}(\omega) \psi^{i}(\d t,\d \omega), \quad x \in \left( \mathcal{L}^{\infty} \right)^{4},
\end{equation}
with each of the $\left( \psi^{i} \right)_{i}$ being $\left( \d t \otimes d\P \right)$-absolutely continuous such that
\begin{equation}
    \left( \d t \otimes d\P \right)(A) = 0 \implies \psi^{i}(A) = 0, \quad A \in \mathcal{B} \left( [0,T] \right) \otimes \mathcal{F}, \quad i \in \{1,2,3,4\}.
\end{equation}
By non-negativity of $l$ on the non-negative ordering cone $\mathcal{C}$ from \eqref{eq:non_negative_separating_linear_form}, it is clear that $\left( \psi^{i} \right)_{i} \in \textbf{\textit{ba}}_{+}^{4}$, and the countably additive and finitely additive parts from their respective Yosida-Hewitt decompositions \eqref{eq:hewitt_yosida_decomposition} are non-negative such that $\left( \lambda^{i} \right)_{i} \in \left( \mathcal{L}_{+}^{1} \right)^{4}$ and $\left( \psi^{i,o} \right)_{i} \in \left( \mathcal{S}_{+} \right)^{4}$ as a consequence of \cite[Theorem 1.23]{yosida1952finitely}. Since $\mathcal{G}(\hat{u}) \in -\mathcal{C}$, equation \eqref{eq:complementary_slackness_functional} readily yields by linearity the complementary slackness conditions \eqref{eq:KKT_slackness_conditions} and \eqref{eq:KKT_slackness_conditions_singular}.\\

\noindent Furthermore, using \eqref{eq:nabla_J}, the stationary equation \eqref{eq:stationary_lagrange} rewrites equivalently as
\begin{align} \label{eq:stationary_intermediate_general_case}
    & \int_{\Omega \times [0,T]} \left( - \alpha_{t}(\omega) + \hat{u}_{t}(\omega) + \int_{0}^{t} K(t,s) \hat{u}_{s}(\omega) + \int_{t}^{T} K(s,t) \E_{t}\hat{u}_{s}(\omega) \d s \right) h_{t}(\omega) \d t \P(\d \omega) \\
    & = \int_{\Omega \times [0,T]} h_{t}(\omega) \left( \psi^{1}(\d t,\d \omega) - \psi^{2}(\d t,\d \omega) \right) + \int_{\Omega \times [0,T]} I(h)(t, \omega) \left( \psi^{3}(\d t,\d \omega) -  \psi^{4}(\d t,\d \omega) \right), \quad  h \in \mathcal{L}^{\infty},
\end{align}
where $I$ is the Lebesgue integral-operator defined by
\begin{equation} \label{eq:definition_lebesgue_integral_operator}
    I(h)(t, \omega) := \int_{0}^{t}h_{s}(\omega)\d s, \quad (t,\omega) \in [0,T] \times \Omega, \quad h \in \mathcal{L}^{\infty}.
\end{equation}
Since $I$ is linear and bounded, we define its unique bounded and linear adjoint operator $I^{\star}: \left( \mathcal{L}^{\infty} \right)^{\star} \to \left( \mathcal{L}^{\infty} \right)^{\star}$ via the dual brackets by
\begin{equation} \label{eq:definition_adjoint_operator}
    \langle I (h), \psi \rangle = \langle h, I^{\star}(\psi) \rangle, \quad  h \in \mathcal{L}^{\infty}, \quad  \psi \in \left( \mathcal{L}^{\infty} \right)^{\star},
\end{equation}
with
\begin{equation}
    \left|\left|I^{*}\right|\right| = \left|\left|I\right|\right|,
\end{equation}
see \cite[Chapter 6, Section 2]{dunford1988linear}. Consequently, equation \eqref{eq:stationary_intermediate_general_case} can be expressed as
\begin{align} \label{eq:stationary_intermediate_general_case_bis}
    & \int_{\Omega \times [0,T]} \left( -\alpha_{t}(\omega) + \hat{u}_{t}(\omega) + \int_{0}^{t} K(t,s) \hat{u}_{s}(\omega) ds + \int_{t}^{T} K(s,t) \E_{t}\hat{u}_{s}(\omega) \d s \right) h_{t}(\omega) \d t \P(\d \omega) \\
    & = \int_{\Omega \times [0,T]} h_{t}(\omega) \left( \psi^{1}(\d t,\d \omega) - \psi^{2}(\d t,\d \omega) + I^{*}\left(\psi^{3}-\psi^{4}\right)(\d t,\d \omega) \right), \quad  h \in \mathcal{L}^{\infty}.
\end{align}
Notice now that using the linearity of $I^{*}$ onto the Yosida-Hewitt decomposition of $\psi^{3}-\psi^{4}$ and applying the Fubini Theorem to the countably additive part with respect to the Lebesgue measure yields
\begin{equation} \label{eq:fubini_on_ba_element}
    \langle h, I^{\star}(\psi^{3}-\psi^{4}) \rangle = \int_{\Omega \times [0,T]} h_{t}(\omega) \left( \int_{t}^{T} \left( \lambda_{s}^{3}(\omega)-\lambda_{s}^{4}(\omega) \right) \d s \right) \d t \mathbb{P}(d \omega) + \langle h, I^{\star}(\psi^{o,3}-\psi^{o,4}) \rangle, \quad h \in \mathcal{L}^{\infty},
\end{equation}
while Fubini Theorem loses its validity in general for the singular part $\psi^{o,3}-\psi^{o,4}$, see \cite[Theorem 3.3]{yosida1952finitely}. Note also that the linearity of $I^{*}$ readily yields the inclusion
\begin{equation} \label{eq:S_stable_by_I_star}
    I^{*} \left( \mathcal{S} \right) \subset \mathcal{S},
\end{equation}
where recall $\mathcal{S}$ is defined just after the Yosida-Hewitt decomposition \eqref{eq:hewitt_yosida_decomposition}. Now, injecting \eqref{eq:fubini_on_ba_element} into \eqref{eq:stationary_intermediate_general_case_bis} and using the tower property yields
\begin{align} \label{eq:stationary_intermediate_general_case_bis_bis}
    & \int_{\Omega \times [0,T]} \left( -\alpha_{t}(\omega) + \hat{u}_{t}(\omega) + \int_{0}^{t} K(t,s) \hat{u}_{s}(\omega) \d s + \int_{t}^{T} K(s,t) \E_{t}\hat{u}_{s}(\omega) \d s \right) h_{t}(\omega) \d t \P(\d \omega) \\
    = & \int_{\Omega \times [0,T]} h_{t}(\omega) \left( \lambda_{s}^{1}(\omega) - \lambda_{s}^{2}(\omega) + \mathbb{E}_{.} \left[ \int_{.}^{T} \left( \lambda_{s}^{3}(\omega) - \lambda_{s}^{4}(\omega) \right) \d s \right] \right) \d t \P(\d \omega) \\
    & + \int_{\Omega \times [0,T]} h_{t}(\omega) \left( \psi^{1,o}(\d t,\d \omega) - \psi^{2,o}(\d t,\d \omega) + I^{*}\left(\psi^{3,o}-\psi^{4,o}\right)(\d t,\d \omega)\right), \quad  h \in \mathcal{L}^{\infty}.
\end{align}
Given that $\alpha \in \mathcal{L}^{2}$ (recall the drift-signal $P$ in \eqref{eq:alpha_signal_definition} is assumed to lie in $\mathcal{L}^{2}$), and since $K$ satisfies \eqref{eq:square_integrable_kernel_assumption}, then
\begin{equation}
   -\alpha_{t}(\omega) + \hat{u}_{t}(\omega) + \int_{0}^{t} K(t,s) \hat{u}_{s}(\omega) \d s + \int_{t}^{T} K(s,t) \E_{t}\hat{u}_{s}(\omega) \d s \in \mathcal{L}^{2},
\end{equation}
and \eqref{eq:stationary_intermediate_general_case_bis_bis} being true for any $\mathcal{L}^{\infty}$-test function $h$, then there is an \textit{aggregated Lagrange multiplier} $\varphi \in \mathcal{L}^{2}$ such that, by uniqueness of the Yosida-Hewitt decomposition, equations \eqref{eq:aggregated_lagrange_multiplier} and \eqref{eq:aggregated_lagrange_multiplier_singular} are satisfied. Consequently, combining \eqref{eq:stationary_intermediate_general_case_bis_bis} and \eqref{eq:aggregated_lagrange_multiplier}, we obtain the desired stochastic Fredholm equation \eqref{eq:Fredholm_equation}.
\end{proof}

\subsection{Putting everything together to prove Theorem~\ref{T:generalized_kkt_theorem_necessary}} \label{ss:proof_necessary_condition}

We are now ready to complete the proof of Theorem~\ref{T:generalized_kkt_theorem_necessary}.
\begin{proof}[Proof of Theorem \ref{T:generalized_kkt_theorem_necessary}] An application of 
Theorem~\ref{T:existenceuniqueness} yields the existence and uniqueness of an admissible optimal control $\hat u$. Assuming in addition that $int(\mathcal U)\neq \emptyset$, Lemma~\ref{L:lagrangemultiplierandFredholm} yields the  existence of Lagrange multipliers $\lambda := \left(\lambda^{i}\right)_{i \in \{1,2,3,4\}} \in \left( \mathcal{L}_{+}^{1} \right)^{4}$ satisfying the decomposition \eqref{eq:aggregated_lagrange_multiplier} with $\varphi \in \mathcal{L}^{2}$, as well as the complementary slackness equations \eqref{eq:KKT_slackness_conditions}, and such that the optimal control $\hat u$  solves the Fredholm equation 
\eqref{eq:Fredholm_equation}. 
Given $K \in \mathcal K$, then such Fredholm equation indeed admits the unique solution $\hat{u}$ explicitly given by \eqref{eq:explicit_optimal_control} as a direct application of \cite[Theorem 2.8]{abijaberOptPort2024} in the uni-dimensional case.
\end{proof}

\subsection{Proof of Theorem \ref{T:generalized_kkt_theorem_sufficient}} \label{ss:proof_sufficient_condition}

\begin{proof}[Proof of Theorem \ref{T:generalized_kkt_theorem_sufficient}]
Fix $\lambda := \left( \lambda^{i} \right)_{i \in \{1,2,3,4\}} \in \left(\mathcal{L}_{+}^{1}\right)^{4}$ such that the process $\varphi$ given by \eqref{eq:aggregated_lagrange_multiplier} lies in $\mathcal{L}^{2}$. Applying again \cite[Theorem 2.8]{abijaberOptPort2024} in the uni-dimensional case, we fix $u^{*} \in \mathcal{L}^{2}$ the unique solution to the Fredholm equation \eqref{eq:assume_fredholm}. Assume now that $u^{*} \in \mathcal U$ and that the couple $(u^{*},\lambda)$ satisfies the complementary slackness conditions \eqref{eq:KKT_slackness_conditions}. We argue that
\begin{equation}\label{eq:tempsuccesive}
    \mathcal{J}(u^{*}) = \mathcal{L}(u^{*}, \lambda) \leq \mathcal{L}(u, \lambda) \leq \mathcal{J}(u), \quad u \in \mathcal{U},
\end{equation}
where $\mathcal{L}$ is the Lagrangian introduced in \eqref{eq:def_lagrange_functional}. Indeed, the first equality holds from the complementary slackness equations \eqref{eq:KKT_slackness_conditions} which can be re-written as $\langle \mathcal{G}(u^{*}), \lambda \rangle = 0$. The first inequality is obtained by noting that $u^*$ is a minimizer of $\mathcal{L}(\cdot, \lambda)$. To see this, we observe that  $\nabla_{u} \mathcal{L}(u^{*}, \psi) = 0$ since $u^{*}$ solves the Fredholm equation \eqref{eq:assume_fredholm}, recall \eqref{eq:nabla_J}-\eqref{eq:nabla_G}. Combined with the fact that  $u \to \mathcal{L}(u, \lambda)$ is convex (recall that $\mathcal{J}$ is $1$-convex from point $(iii)$ in Lemma~\ref{L:properties_functional_admin_controls}, and that $u\to \langle u,\lambda \rangle_{\mathcal{Y}, \left( \left( \mathcal{L}^{\infty} \right)^{4} \right)^{\star}}$ is linear in $u$), this yields that $\mathcal{L}(\cdot, \lambda)$ reaches a global minimum at $u^{*}$.  Finally, the second inequality in \eqref{eq:tempsuccesive} is a consequence of the fact that
\begin{equation}
    \langle \mathcal{G}(u), \lambda \rangle \leq 0, \quad u \in \mathcal{U}.
\end{equation}
since $\lambda^{i} \geq 0, \; i \in \{1,2,3,4 \}$ and  $\mathcal{G}(u) \in -\mathcal{C}, \; u \in \mathcal{U}$.

Since $u^*\in \mathcal U$ by assumption, it follows from \eqref{eq:tempsuccesive} that $u^*$ is a minimizer of $\mathcal J$ over $\mathcal U$. Then, by uniqueness of the minimizer of $\mathcal{J}$ on $\mathcal{U}$ from Theorem~\ref{T:existenceuniqueness}, we obtain that $u^*=\hat u$, which concludes the proof.
\end{proof}

\end{document}